\documentclass[12pt,twoside,english,a4paper]{article}
\usepackage{etex}

\usepackage{caption}
\usepackage[latin1]{inputenc}
\usepackage[intlimits] {amsmath}   
\usepackage{graphicx}
\usepackage{psfrag}
\usepackage{ifthen}
\usepackage{fancyhdr}
\usepackage{rotating}
\usepackage{multirow}
\usepackage{booktabs}              
\usepackage{amssymb}
\usepackage{amsbsy}
\usepackage{bm}                    
\usepackage{bbm}
\usepackage{babel}
\usepackage{theorem}
\usepackage{upgreek}  
\usepackage{pifont}   
\usepackage[active]{srcltx}   

\usepackage{suetterl}  
\newfont{\suetdbl}{suet14 scaled 2000}  

\usepackage{yfonts}  
\newfont{\gothdbl}{ygoth scaled 2000} 
\newfont{\frakdbl}{yfrak scaled 2000}
\newfont{\swabdbl}{yswab scaled 2000}

\usepackage[normalem]{ulem}  
\usepackage[textsize=footnotesize,textwidth=1.7cm]{todonotes}
\usepackage[]{defjs5}      

\usepackage{float}
\usepackage{caption}
\usepackage{wrapfig}
\usepackage{url}

\usepackage{pgfplots}

\pgfplotsset{compat=1.18}

\usepackage[authoryear]{natbib}
\bibpunct{[}{]}{,}{n}{}{;}
\fancypagestyle{plain}{%
  \fancyhead{}%
  \fancyfoot[c]{\sffamily\thepage}%
}
\makeatletter                           
\def\cleardoublepage{\clearpage\if@twoside \ifodd\c@page\else
  \hbox{}
  \vspace*{\fill}
  \thispagestyle{empty}
  \newpage
  \if@twocolumn\hbox{}\newpage\fi\fi\fi}
\makeatother
\begin{document}
\unitlength1.0cm
\frenchspacing


\newpage

\thispagestyle{empty}

\ce{\bf \large
A third-medium approach for electro-thermo-mechanical contact
}
\vspace{2mm}
\ce{\bf \large
considering Joule heating
}

\vspace{4mm}
\ce{M. Vorwerk}

\vspace{4mm}
\ce{Institute of Mechanics, Faculty of Engineering, University of Duisburg-Essen}
\ce{\small e-mail: maximilian.vorwerk@uni-due.de}

\vspace{4mm}
\begin{center}
{\bf \large Abstract}
\bigskip

{\footnotesize
\begin{minipage}{16.cm}
\noindent
Connectors are ubiquitous in technical systems and frequently combine mechanical contact with electric current transfer and heat generation.
To capture these interacting processes, an electro-thermo-mechanical contact formulation based on the third-medium concept is proposed.
Mechanical contact, electric current flow, and heat conduction are represented within a unified finite element framework without explicit contact-surface tracking.
Deformation-dependent switching functions govern the onset of electrical and thermal transport during interface closure.
Electrical conduction is coupled to the transient thermal problem through Joule heating, while the temperature dependence of the electrical conductivity provides the corresponding feedback on the current flow.
The coupled displacement, electric potential, and temperature fields are solved monolithically and combined with an independently interpolated deformation-gradient-like field for robust third-medium regularization under severe compression.
Numerical examples demonstrate contact-induced current transfer, subsequent Joule heating and transient temperature evolution, as well as localized current paths during progressive closure of rough interfaces.
\end{minipage}
}
\end{center}


\section{Introduction}
\vspace{-4mm}

The finite element method provides a central tool for simulation-driven engineering design.
Contact problems remain particularly demanding since contact zones, boundary conditions and localized force transfer evolve during deformation, cf.~\cite{Wri:2006:ccm,Lau:2002:cca,Yas:2013:nmi}.
Relevant applications range from metal forming and crash analysis to soft robotics and contact-aided mechanisms, cf.~\cite{RusTol:2015:dfa,FreSigPou:2024:too}.
Classical formulations enforce contact constraints on explicitly discretized interfaces using penalty, Lagrange multiplier, augmented Lagrangian or barrier methods, cf.~\cite{Wri:2006:ccm,Lau:2002:cca}.
Surface-to-surface and mortar formulations extend this framework to non-matching meshes, cf.~\cite{PusLau:2004:ams,SauDeL:2013:acc}.
Contact search, gap evaluation and active-set treatment become increasingly involved for large sliding, self-contact and topology optimization with initially unknown contact boundaries, cf.~\cite{BluSigPou:2021:icm,FreSigPou:2024:too,WriKorJun:2025:atm}.\\
Additional challenges arise when mechanical contact is coupled to electric current and heat transfer.
Early electro-thermo-mechanical formulations account for contact-dependent electrical resistance, Joule heating and the resulting temperature evolution, cf.~\cite{WeiWri:2010:nme}.
Pressure-dependent electrical and thermal transfer laws have also been developed for electrically assisted forming processes, cf.~\cite{TerOzhTraShiMatSolKlo:2017:etm}.
Fully coupled finite element models resolve electrical contact resistance and the corresponding Joule heating, cf.~\cite{LiSheKe:2022:mpe}.
Rough-surface effects and their influence on current transfer and temperature rise are considered in \cite{LiSheGueKe:2024:aem}.
Most available approaches rely on explicitly defined contact interfaces and retain the classical requirements of contact detection, gap evaluation and active-set treatment.\\
Third-medium contact offers an alternative by replacing explicit interface constraints with a highly compliant fictitious material occupying the space between potentially contacting bodies.
Its negligible stiffness limits the mechanical influence before contact, whereas strong compression enables the transmission of contact forces without contact search or inequality constraints.
An initial finite-deformation formulation is presented in \cite{WriSchSch:2013:afe}.
Extensions to high-order finite elements, isogeometric analysis and isogeometric-meshfree coupling are developed in \cite{BogZanKolRan:2015:ncw,KruNguWriDeL:2018:ifc,HuaNguZho:2018:aim}.
Particular benefits arise in density-based topology optimization, where the compliant void phase can simultaneously act as a contact medium.
Internal contact, self-contacting structures and tailored nonlinear springs are addressed in \cite{BluSigPou:2021:icm,FreSigPou:2024:too,BluSigPou:2023:ido}.
Thermo-mechanical regulators and frictional contact extend the range of accessible applications, cf.~\cite{DalAleFrePouSig:2025:too,FreRokPouSigGee:2024:aft}.
Three-dimensional hyperelastic contact and pneumatically actuated systems are considered in \cite{XuXueWri:2026:tdt}.
Robustness under extreme compression remains governed by the stabilization of the fictitious medium.
The HuHu regularization controls severe distortions through displacement Hessians, cf.~\cite{BluSigPou:2021:icm}.
Its HuHu-LuLu extension reduces excessive penalization of bending and quadratic compression, cf.~\cite{FreDalSigPou:2025:itm}.
Alternative strategies employ deformation-gradient averaging or rotation-based measures, cf.~\cite{FalAmaHor:2026:dga,DahSjoDalWal:2026:ara}.
Recent developments increasingly target low-order discretizations.
Auxiliary-field formulations for first-order finite elements are proposed in \cite{WriKorJun:2025:atm}.
Thermo-mechanical transport through the third-medium is considered in \cite{Wri:2026:atm}.
Virtual-element formulations permit polygonal meshes without conventional stabilization terms, cf.~\cite{XuWri:2026:sfv}.
A first-order virtual-element extension is presented in \cite{XuXueWri:2026:fov}.
A neighbored-element strategy reduces the number of globally coupled regularization unknowns, cf.~\cite{ZabJanJun:2026:afa}.
An independently interpolated deformation-gradient-like field provides another low-order-compatible stabilization and avoids the direct evaluation of second displacement derivatives, cf.~\cite{VorSchWri:2026:arm}.

Despite these advances, coupled electric current transfer, Joule heating and transient heat conduction across closing interfaces have received little attention within third-medium contact.
The central contribution of the present work is the extension of the third-medium concept to electrically conducting contact with Joule heating and temperature-dependent electrical conductivity.
Electrical contact formation is represented through a deformation-dependent activation law, while transient heat conduction accounts for the subsequent thermal evolution.
The formulation is embedded into the mixed low-order third-medium framework of~\cite{VorSchWri:2026:arm}, which provides the mechanical regularization required under severe compression.
Numerical examples assess contact-induced current transfer, the onset of Joule heating and the resulting transient thermal response.

The remainder of the paper introduces the coupled continuum formulation and the constitutive response of the conducting solids and the third-medium.
Subsequent sections present the switching functions, finite element discretization and numerical examples, followed by the main conclusions.

\section{Continuum formulation of deformable conductors}
\vspace{-4mm}
Within this section, the thermo-electro-mechanical continuum formulation of a deformable conductor is introduced. 
The formulation combines quasi-static mechanical equilibrium, stationary electric current conduction and transient heat transport under finite deformations. 
Its subsequent extension to the conducting third-medium is introduced separately.
The physical solid domain is denoted by $\B_{\mathrm{s}}$, while the space between potentially contacting surfaces is represented by the fictitious third-medium domain $\B_{\mathrm{tm}}$. 
Both domains are described within finite deformation kinematics. 
For the solid, a standard hyperelastic material model is employed. 
In the third-medium, the same hyperelastic base energy is used in a strongly scaled form and supplemented by an additional stabilization contribution.  

\subsection{Boundary value problem, kinematics and weak form of a conductor}
The coupled boundary value problem is formulated in the reference configuration $\B\subset\IR^3$ and combines quasi-static mechanical equilibrium, stationary electric current conservation and transient heat transport as
\begin{equation}
	\Div{\bP}
	+
	\overline \bb
	=
	\bzero,
	\quad
	\Div{\bJ}_{\mathrm{e}}
	=
	0,
	\quad\textrm{and}\quad	
	\rho_0 c_0 \dot{\theta}
	+
	\Div{\bQ}-Q_J
	=
	0
	\quad\textrm{in}\quad\B,
	\label{eq:FieldEq}
\end{equation}
where these equations denote the mechanical equilibrium equation, the stationary current balance and the transient heat balance, respectively. 
Here, $\bP$ denotes the first Piola-Kirchhoff stress tensor, $\rho_0$ is the referential mass density, $\overline{\bb}$ denotes the body force per unit reference volume, $\bJ_{\mathrm{e}}$ is the referential electric current density, $\bQ$ denotes the referential heat flux, and $Q_J$ is the volumetric Joule heat source per unit reference volume. 
The parameter $c_0$ is the specific heat capacity. 
To complete the definition of the boundary value problem, the boundary conditions are prescribed as
\begin{align}
	\bu
	&=
	\bar{\bu}
	&& \mathrm{on}\quad \partial\B_u ,
	&
	\bP\cdot\bN
	&=
	\bar{\bt}
	&& \mathrm{on}\quad \partial\B_t ,
	\label{eq:BC_mech}
	\\
	\varphi
	&=
	\bar{\varphi}
	&& \mathrm{on}\quad \partial\B_\varphi ,
	&
	\bJ_{\mathrm{e}}\cdot\bN
	&=
	\bar{j}_N
	&& \mathrm{on}\quad \partial\B_{\bar{j}_N} ,
    \\
	\theta
	&=
	\bar{\theta}
	&& \mathrm{on}\quad \partial\B_\theta ,
	&
	\bQ\cdot\bN
	&=
	\bar{Q}_N
	&& \mathrm{on}\quad \partial\B_{\bar{Q}_N} .    	
	\label{eq:BC}
\end{align}
Here, $\bar{\bu}$, $\bar\varphi$ and $\bar{\theta}$ are the prescribed displacement, electric potential and temperature, while $\bar{\bt}$, $\bar{j}_N$ and $\bar{Q}_N$ are the prescribed traction, normal current flux and normal heat flux.
$\bN$ denotes the outward unit normal in the reference configuration.
The transient heat equation is supplemented by the initial condition  $\theta(\bX,0)=\theta_{\mathrm{init}}(\bX)$ in $\B$.
To obtain the weak formulation, the field equations in Eq.~\ref{eq:FieldEq} are multiplied by the admissible test functions $\delta\bu$, $\delta\varphi$ and $\delta\theta$ and integrated over the reference configuration. 
Application of integration by parts yields the coupled variational problem
\begin{equation}
\begin{aligned}
    G(\bu,\varphi,\theta;\delta\bu,\delta\varphi,\delta\theta)=&\,
    G_u(\bu,\theta;\delta\bu)+
    G_\varphi(\bu,\varphi,\theta;\delta\varphi)+
    G_\theta(\bu,\varphi,\theta;\delta\theta)=0\\[3mm]
	G_u(\bu,\theta;\delta\bu)
	=&
	\int_{\B}
	\bP
	:
	\delta\bF
	\,\mathrm{d}V
	-
	\int_{\B}
	\overline
	\bb
	\cdot
	\delta\bu
	\,\mathrm{d}V
	-
	\int_{\partial\B_t}
	\bar{\bt}
	\cdot
	\delta\bu
	\,\mathrm{d}A\\[1mm]
	G_\varphi(\bu,\varphi,\theta;\delta\varphi)
	=&
	\int_{\B}
	\bJ_e
	\cdot
	\delta\bE	
	\,\mathrm{d}V
	+
	\int_{\partial\B_{\bar{j}_N}}
	\delta\varphi\,
	\bar{j}_N
	\,\mathrm{d}A\\[1mm]
	G_\theta(\bu,\varphi,\theta;\delta\theta)
	=
	&
	-
	\int_{\B}
	\nabla_X\delta\theta
	\cdot
	\bQ
	\,\mathrm{d}V
	+
	\int_{\B}
	\delta\theta\,
	\rho_0c_0\dot\theta
	\,\mathrm{d}V...
	\\
	&\hspace*{30mm}...
	-
	\int_{\B}
	\delta\theta\,
	Q_J
	\,\mathrm{d}V
	+
	\int_{\partial\B_{\bar{Q}_N}}
	\delta\theta\,
	\bar Q_N
	\,\mathrm{d}A.		
\end{aligned}
\label{eq:weak_form}
\end{equation}
Here, $\delta\bF=\nabla_X\delta\bu$, $\delta\bE=-\nabla_X\delta\varphi$ and $\nabla_X\delta\theta$ denote the virtual deformation gradient, the virtual referential electric field and the gradient of the virtual temperature, respectively.
Body forces are neglected throughout the following developments.
A material point $\bX\in\B$ is mapped to the current configuration through $\bx=\Bchi(\bX,t)=\bX+\bu(\bX,t)$ such that the deformation gradient, its Jacobian and the right Cauchy-Green tensor are defined by
\begin{equation}
	\bF
	=
	\frac{\partial \bx}{\partial \bX}
	=
	\bI+\nabla_X\bu ,
	\qquad
	J
	=
	\det\bF
	\qquad\textrm{and}\qquad
	\bC
	=
	\bF^{T}\cdot\bF .
	\label{eq:F}
\end{equation}
The isochoric part of the deformation is represented by the volume-preserving deformation gradient and the corresponding right Cauchy-Green tensor,
\begin{equation}
	\widehat{\bF}
	=
	J^{-1/3}\bF
	\qquad\textrm{and}\qquad
	\widehat{\bC}
	=
	\widehat{\bF}^{T}\cdot\widehat{\bF}
	=
	J^{-2/3}\bC .
\end{equation}
Thermal expansion is incorporated through a multiplicative decomposition of the deformation gradient into elastic and thermal parts,
\begin{equation}
	\bF
	=
	\bF_e
	\cdot
	\bF_\theta ,
	\label{eq:Fsplit}
\end{equation}
where $\bF_e$ and $\bF_\theta$ denote the elastic and thermal deformation gradients, respectively.
Transient thermal response follows from the balance of energy, where the heat flux is described by Fourier's law.
In the present formulation, a separate thermal free-energy contribution is not introduced, since the transient heat equation is formulated directly in terms of the specific heat capacity and Fourier heat flux.
Assuming isotropic thermal expansion, the thermal deformation gradient is defined as
\begin{equation}
	\bF_{\theta}
	=
	\lambda_{\theta}\bI,
	\qquad
	\lambda_{\theta}
	=
	1+\alpha_{\theta}(\theta-\theta_0),
\end{equation}
where $\alpha_{\theta}$ denotes the coefficient of thermal expansion and $\theta_0$ the reference temperature.
For the plane-strain problems considered in this work, thermal expansion is restricted to the in-plane directions, such that $\bF_{\theta}=\operatorname{diag}(\lambda_{\theta},\lambda_{\theta},1)$.
The linear relation is adopted for moderate temperature changes around the reference state.
The elastic deformation gradient and the corresponding elastic right Cauchy-Green tensor are therefore obtained as
\begin{equation}
	\bF_e
	=
	\bF
	\cdot
	\bF_\theta^{-1}
\qquad\textrm{and}\qquad
	\bC_e
	=
	\bF_e^{T}
	\cdot
	\bF_e.
	\label{eq:FeandCe}
\end{equation}
Mechanical constitutive relations are formulated in terms of $\bC_e$. 
Consequently, thermal expansion affects the stress response through the elastic deformation gradient.

\subsection{Behavior of the conducting contact domains}
The conducting solids are modeled as deformable metallic conductors.
Deformation-induced changes in the electrical conductivity are assumed to be negligible, while Joule heat is transported through the conductor by Fourier heat conduction.
The mechanical response is described by the Helmholtz free-energy density $\psi_M$.
Electrical conduction is introduced separately through a dissipation potential.
Heat conduction, Joule heating and heat storage enter through the transient energy balance.
No additional thermal contribution is therefore included in the Helmholtz free-energy density.

\textbf{The mechanical response}
is described by a compressible Neo-Hookean material model.
In the conducting solid, the strain-energy density is given by
\begin{equation}
	\psi_M^s(\bC_e)
	=
	\frac{K}{2}
	\left(
	\ln J_e
	\right)^2
	+
	\frac{\mu}{2}
	\left(
	J_e^{-2/3}
	\tr\bC_e
	-
	3
	\right)
	\qquad\mathrm{with}\qquad
	J_e
	=
	\det\bF_e,
	\label{eq:NeoHook}
\end{equation}
where $K$ and $\mu$ denote the bulk and shear moduli, respectively.
Dependence on the elastic right Cauchy-Green tensor $\bC_e$ and the elastic Jacobian $J_e$ accounts for the mechanically effective part of the deformation.
The first Piola-Kirchhoff stress follows by the chain rule from the dependence of the free energy on $\bC_e$, yielding
\begin{equation}
	\bS
	=
	2
	\frac{\partial\psi_M^s}{\partial\bC_e}
	\qquad\mathrm{and}\qquad
	\bP
	=
	\bF_e
	\cdot
	\bS
	\cdot
	\bF_\theta^{-T},
	\label{eq:PK}
\end{equation}
where $\bS$ and $\bP$ denote the second and first Piola-Kirchhoff stress tensors, respectively.

\textbf{The electric response}
is introduced through the electrical dissipation potential
\begin{equation}
	\mathcal{D}_{\varphi}^{s}(\bE)
	=
	\frac{1}{2}
	\bE
	\cdot
	\bK_e^{\,s}
	\cdot
	\bE ,
	\label{eq:ElectricalEnergy}
\end{equation}
where $\bK_e^{\,s}$ denotes the effective electrical conductivity tensor of the solid in the reference configuration.
The spatial electric field and current density follow from the electric potential as
\begin{equation}
	\be
	=
	-
	\nabla_x\varphi
	\qquad\textrm{and}\qquad
	\bj_{\mathrm e}
	=
	\sigma_{\mathrm s}
	\be ,
	\label{eq:EandDensity}
\end{equation}
where $\sigma_{\mathrm s}$ denotes the spatial electrical conductivity of the conducting solid.
Piezoresistive changes in $\sigma_{\mathrm s}$ are neglected, while the geometric influence of the deformation on electric transport is retained.
Application of the Piola transformation yields
\begin{equation}
	\bJ_{\mathrm e}^{\,s}
	=
	J
	\bF^{-1}
	\cdot
	\bj_{\mathrm e}
	=
	J
	\sigma_{\mathrm s}
	\bC^{-1}
	\cdot
	\bE
	=
	\bK_e^{\,s}
	\cdot
	\bE
	\qquad\textrm{with}\qquad
	\bK_e^{\,s}
	=
	J
	\sigma_{\mathrm s}
	\bC^{-1},
	\label{eq:OhmLawReferential}
\end{equation}
where $\bJ_{\mathrm e}^{\,s}$ denotes the referential current density and
$\bE=-\nabla_X\varphi$ is the referential electric field.
Although the scalar conductivity $\sigma_{\mathrm s}$ remains independent of the elastic strain, the referential conductivity tensor depends on the deformation through $J$ and $\bC^{-1}$.
The constitutive current relation follows from the dissipation potential as
\begin{equation}
	\bJ_{\mathrm e}^{\,s}
	=
	\frac{
		\partial
		\mathcal{D}_{\varphi}^{s}
	}{
		\partial
		\bE
	}.
	\label{eq:SimplifiedConductingMaterial}
\end{equation}
The corresponding extension to the conducting third-medium is introduced in the following subsection.

\textbf{The thermal response}
is governed by the transient heat balance introduced in Eq.~\ref{eq:FieldEq}.
Heat storage is described directly by the constant specific heat capacity $c_0$, while heat transport follows Fourier's law.
The thermal dissipation potential in the reference configuration is defined as
\begin{equation}
	\mathcal{D}_{\theta}^{s}
	\left(
		\nabla_X\theta
	\right)
	=
	\frac{1}{2}
	\nabla_X\theta
	\cdot
	\bK_\theta^{\,s}
	\cdot
	\nabla_X\theta ,
	\label{eq:ThermalDissipationPotential}
\end{equation}
where $\bK_\theta^{\,s}$ denotes the effective thermal conductivity tensor of the solid in the reference configuration.
The spatial heat flux follows as
\begin{equation}
	\bq
	=
	-
	k_{\mathrm s}
	\nabla_x\theta ,
	\label{eq:FourierSpatial}
\end{equation}
where $k_{\mathrm s}$ denotes the spatial thermal conductivity of the conducting solid.
Application of the Piola transformation yields
\begin{equation}
	\bQ^{s}
	=
	J
	\bF^{-1}
	\cdot
	\bq
	=
	-
	J
	k_{\mathrm s}
	\bC^{-1}
	\cdot
	\nabla_X\theta
	=
	-
	\bK_\theta^{\,s}
	\cdot
	\nabla_X\theta
	\qquad\mathrm{with}\qquad
	\bK_\theta^{\,s}
	=
	J
	k_{\mathrm s}
	\bC^{-1}.
	\label{eq:FourierReferential}
\end{equation}
The referential heat flux follows from the thermal dissipation potential according to
\begin{equation}
	\bQ^{s}
	=
	-
	\frac{
		\partial
		\mathcal{D}_{\theta}^{s}
	}{
		\partial
		\nabla_X\theta
	}.
	\label{eq:ThermalDissipationRelation}
\end{equation}
Electric current generates heat through Joule dissipation with the corresponding referential volumetric heat source given by
\begin{equation}
	Q_J^{s}
	=
	\bJ_{\mathrm e}^{\,s}
	\cdot
	\bE
	=
	\bE
	\cdot
	\bK_e^{\,s}
	\cdot
	\bE
	=
	2
	\mathcal{D}_{\varphi}^{s}.
	\label{eq:JouleHeat}
\end{equation}
Accordingly, $Q_J^{\mathrm{s}}\geq 0$ denotes the positive volumetric heat-generation rate associated with electrical dissipation.
This generated heat is subsequently redistributed within the conducting solid by Fourier heat conduction.
Heat conduction is therefore treated as an irreversible transport process, while heat storage enters directly through the volumetric heat-capacity term $\rho_0 c_0 \dot{\theta}$ in the transient heat balance.
The energetic and dissipative contributions may be collected in a generalized variational
functional, but their sum must not be interpreted as a stored Helmholtz free energy.

\subsection{Behavior of the current and temperature conducting third-medium}
The third-medium differs fundamentally from the conducting solids, since it represents a fictitious contact material rather than a physical conductor.
Besides providing the mechanical contact barrier, it enables electrical and thermal transport only after sufficient compression has occurred.
Mechanical, electrical and thermal constitutive relations are therefore modified to account for the evolving contact state.

\textbf{The third-medium mechanical response}
is based on the same hyperelastic constitutive model introduced for the conducting solids.
Since the third-medium represents a fictitious contact material rather than a physical phase, only a negligible mechanical stiffness is assigned before contact.
Accordingly, the strain-energy density is scaled by the small parameter $\gamma$ to
\begin{equation}
	\psi_M^{\,tm}
	=
	\gamma\,
	\psi_M .
	\label{eq:TMCHyperelastic}
\end{equation}
Compression of the third-medium therefore produces only a negligible mechanical response under moderate deformations, while the rapidly increasing Neo-Hookean energy prevents the complete collapse of the intermediate layer.
For the three-dimensional formulation, the complete Neo-Hookean strain-energy density introduced in Eq.~\eqref{eq:NeoHook} is retained.
Within the present two-dimensional plane-strain setting, only the isochoric contribution is considered (cf.~\cite{WriKorJun:2025:atm}), resulting in
\begin{equation}
	\psi_M^{\,tm}
	=
	\gamma
	\frac{\mu}{2}
	\left(
	J_e^{-2/3}
	\operatorname{tr}
	\bC_e
	-
	3
	\right).
	\label{eq:TMCReducedEnergy}
\end{equation}
This specialization is restricted to the present two-dimensional plane-strain formulation.
Embedding the two-dimensional deformation into three dimensions by setting $F_{33}=1$ ensures that the isochoric contribution alone provides the desired contact-barrier effect, since the energy becomes unbounded for $J\rightarrow0$.
The thermo-mechanical decomposition introduced in Eq.~\ref{eq:Fsplit} is retained within the third-medium, such that thermal expansion also enters its mechanical response through $\bF_{\mathrm e}$ and $\bC_{\mathrm e}$.
Its mechanical influence remains strongly reduced by the scaling parameter $\gamma$.

\textbf{The third-medium mixed stabilization behavior}
is introduced to improve the robustness of the third-medium formulation under severe compressive deformations.
Low-order finite element discretizations may otherwise suffer from excessive element distortions once the intermediate layer approaches complete collapse.
Following the mixed third-medium formulation proposed in~\cite{VorSchWri:2026:arm}, an additional deformation-gradient-like field $\BTheta$ is introduced and interpolated independently of the displacement field.
This auxiliary field approximates the deformation gradient according to $\BTheta\approx\bF$.
Weak coupling between both fields is enforced through the penalty contribution
\begin{equation}
	W_{p}^{\mathrm{tm}}
	=
	\int_{\B_{\mathrm{tm}}}
	p_\Theta
	\left\|
	\BTheta-\bF
	\right\|^2
	\,\mathrm{d}V,
	\label{eq:TMCPenalty}
\end{equation}
where $p_\Theta$ denotes the penalty parameter.
For sufficiently large values of $p_\Theta$, the auxiliary field converges towards the deformation gradient while remaining an independent finite element field.
Additional smoothing is introduced through the gradient regularization
\begin{equation}
	W_{r}^{\mathrm{tm}}
	=
	\int_{\B_{\mathrm{tm}}}
	\alpha_r
	\left\|
	\nabla_X\BTheta
	\right\|^2
	\,\mathrm{d}V,
	\label{eq:TMCRegularization}
\end{equation}
where $\alpha_r$ denotes the regularization parameter.
Since only first derivatives of $\BTheta$ are required, second derivatives of the displacement field are avoided completely.
The stabilized third-medium energy is therefore given by
\begin{equation}
\begin{aligned}
	W_{\mathrm{stab}}^{\mathrm{tm}}
	&=
	W_{M}^{\mathrm{tm}}
	+
	W_{p}^{\mathrm{tm}}
	+
	W_{r}^{\mathrm{tm}}
	\\[2mm]
	&=
	\int_{\B_{\mathrm{tm}}}
	\left[
	\psi_{M}^{\mathrm{tm}}
	+
	p_\Theta
	\left\|
	\BTheta-\bF
	\right\|^2
	+
	\alpha_r
	\left\|
	\nabla_X\BTheta
	\right\|^2
	\right]
	\,\mathrm{d}V .
\end{aligned}
\label{eq:TMC_total_energy_theta}
\end{equation}
The proposed mixed formulation preserves the efficiency of low-order finite elements while substantially improving robustness under extreme compressive deformations of the third-medium.

\textbf{The third-medium electric response}
depends on both the deformation and the temperature of the third-medium.
Compression establishes electrically conducting paths between opposing solid surfaces, whereas separation suppresses current flow.
Accordingly, the constitutive relation is written as
\begin{equation}
	\bJ_{\mathrm e}^{\,tm}
	=
	\bK_e^{\,tm}
	\cdot
	\bE,
	\qquad\textrm{with}\qquad
	\bK_e^{\,tm}
	=
	J
	\sigma_{\,tm}
	\,
	\bC^{-1},
	\label{eq:TMCConductingMaterial}
\end{equation}
where $\sigma_{\,tm}$ denotes the effective electrical conductivity of the third-medium.
Temperature dependence is inherited from the conducting solid, while the deformation-dependent activation of electrical transport is described by
\begin{equation}
\sigma_{\mathrm{tm}}(\theta,J)
=
\sigma_{\mathrm s}(\theta)\,
f_{\varphi}(J),
	\label{eq:sigma_tm_general}
\end{equation}
where $f_\varphi(J)$ denotes the electrical switching function and satisfies
$0\leq f_\varphi(J)\leq1$.
Vanishing values represent an electrically open contact, whereas values approaching unity recover the temperature-dependent conductivity of the conducting solid.
A specific definition of $f_\varphi(J)$ is introduced in the following section.
The electric third-medium dissipation potential follows in analogy to Eq.~\ref{eq:ElectricalEnergy} but involves the constitutive relation and conductivity of the third-medium presented in Eq.~\ref{eq:TMCConductingMaterial} and Eq.~\ref{eq:sigma_tm_general}.

\textbf{The third-medium thermal response}
follows a related concept, but employs a separate switching function.
Heat transfer across the third-medium increases as the distance between opposing solid surfaces decreases.
Accordingly, the referential heat flux is written as
\begin{equation}
	\bQ^{\,tm}
	=
	-
	\bK_\theta^{\,tm}
	\cdot
	\nabla_X\theta,
	\qquad\textrm{with}\qquad
	\bK_\theta^{\,tm}
	=
	J
	k_{\,tm}\,
	\bC^{-1},
	\label{eq:TMCHeatConduction}
\end{equation}
where $\bK_\theta^{\,tm}$ denotes the effective thermal conductivity tensor of the third-medium in the reference configuration.
Effective thermal conductivity is defined by
\begin{equation}
	k_{\,tm}
	=
	k_{\,s}
	\,f_\theta(J),
	\label{eq:TMCHeatConductivity}
\end{equation}
where $k_s$ denotes the thermal conductivity of the conducting solid and $f_\theta(J)$ represents the thermal switching function.
In contrast to $f_\varphi(J)$, the function $f_\theta(J)$ allows thermal transport to evolve according to a separate deformation-dependent activation law.
Both switching functions are specified in detail in the following section.
The thermal third-medium dissipation potential follows in analogy to Eq.~\ref{eq:ThermalDissipationPotential} but involves the constitutive relation and conductivity of the third-medium presented in Eq.~\ref{eq:TMCHeatConduction} and Eq.~\ref{eq:TMCHeatConductivity}.

\subsection{Switching functions within the third-medium \label{sec:SwitFunc}}
Electrical and thermal transport inside the third-medium are activated by two different switching functions.
Although both mechanisms depend on the local compression state, their physical characteristics differ.
Electrical current requires a continuous conducting path between opposing conductors.
Consequently, electrical transport is activated only after sufficient compression of the third-medium.
A discontinuous switching function
\begin{equation}
	f_{\varphi}(J)
	=
	\begin{cases}
	1,
	&
	J<J_{\mathrm{crit}},
	\\[1mm]
	0,
	&
	J\geq J_{\mathrm{crit}} ,
	\end{cases}
	\label{eq:CurrentSwit}
\end{equation}
is therefore adopted to represent the transition from an electrically open to a closed contact.
The discontinuous switching function is treated by an iterative active-set-type update, while the local electrical state of the third-medium is kept fixed within each Newton iteration and updated after completion of the iteration based on the current deformation state. 
Electrical conduction is activated for the subsequent iteration once the local volume ratio satisfies $J<J_{\mathrm{crit}}$. 
Consequently, the consistent linearization applies to the constitutive response within a fixed electrical state, whereas changes of the electrical state are treated between Newton iterations. 
Since the electric problem is stationary and the transient evolution enters only through the comparatively slow thermal response associated with Joule heating, this treatment did not cause numerical instabilities in the numerical examples considered below.
It should be emphasized that $J_{\mathrm{crit}}$ does not represent a universal physical gap or contact distance.
Since the local volume ratio depends on the initial geometry and the deformation of the third-medium, the value of $J_{\mathrm{crit}}$ is problem-dependent and must be selected consistently with the particular third-medium configuration.
Accordingly, the value $J_{\mathrm{crit}}=0.01$ used in the numerical examples below should be regarded as a model parameter controlling the onset of electrical conduction rather than as a universal material parameter.
The same geometric dependence applies to the thermal activation function $f_\theta(J)$.
Hence, both transport laws should be interpreted as third-medium constitutive regularizations parameterized by the local compression state rather than as direct representations of a physical interface gap.
A formulation in terms of a geometrically objective physical gap measure would remove this dependence, but is beyond the scope of the present third-medium model.
Heat transfer is assumed to increase continuously as the intermediate layer becomes thinner.
To prevent negative thermal conductivities, the thermal activation function is defined as
\begin{equation}
	f_{\theta}(J)
	=
	\frac{1}{2}
	\left[
	g(J)
	+
	\sqrt{
		g(J)^2
	+
	\varepsilon_{\theta}^{\,2}
	}
	\right]
	\quad\textrm{with}\quad
	g(J)
	=
	\frac{
	\exp(-\beta J)
	-
	\exp(-\beta)
	}{
	1
	-
	\exp(-\beta)
	}.
	\label{eq:thermal_switch}
\end{equation}

\begin{figure}[H]
    \centering
    \begin{tikzpicture}
        \begin{axis}[
            width=9cm,
            height=4.5cm,
            xmin=0,
            xmax=1,
            ymin=-0.05,
            ymax=1.08,
            xlabel={$J$},
            ylabel={$f(J)$},
            axis lines=left,
            xtick={0.01,0.2,0.4,0.6,0.8,1.0},
            xticklabels={$0.01$,$0.2$,$0.4$,$0.6$,$0.8$,$1.0$},
            ytick={0,0.2,0.4,0.6,0.8,1.0},
            legend style={
                at={(0.97,0.97)},
                anchor=north east,
                draw=none,
                fill=white
            },
            legend cell align={left},
            tick align=outside,
            clip=false
        ]

        \addplot[
            blue,
            very thick
        ]
        coordinates {
            (0,1)
            (0.01,1)
        };
        \addlegendentry{$f_{\varphi}(J)$}

        \addplot[
            blue,
            very thick,
            forget plot
        ]
        coordinates {
            (0.01,0)
            (1,0)
        };

        \addplot[
            only marks,
            mark=o,
            mark size=2.5pt,
            very thick,
            blue,
            fill=white,
            forget plot
        ]
        coordinates {
            (0.01,1)
        };

        \addplot[
            only marks,
            mark=*,
            mark size=2.5pt,
            blue,
            forget plot
        ]
        coordinates {
            (0.01,0)
        };

        \addplot[
            red,
            very thick,
            domain=0:1,
            samples=200
        ]
        {
            0.5*(
                (exp(-5*x)-exp(-5))/(1-exp(-5))
                +
                sqrt(
                    ((exp(-5*x)-exp(-5))/(1-exp(-5)))^2
                    +
                    (1e-4)^2
                )
            )
        };
        \addlegendentry{
            $f_{\theta}(J)$,
            $\beta=5$,
            $\varepsilon_{\theta}=10^{-4}$
        }

        \addplot[
            gray,
            dashed,
            thin,
            forget plot
        ]
        coordinates {
            (0.01,-0.05)
            (0.01,1.05)
        };

        \node[
            anchor=south west,
            font=\small
        ]
        at (axis cs:0.27,0.3)
        {$J_{\mathrm{crit}}=0.01$};

        \end{axis}
    \end{tikzpicture}
    \caption{
        Electrical and thermal switching functions within the third-medium.
        Electrical transport is activated discontinuously for
        $J<J_{\mathrm{crit}}=0.01$, whereas thermal transport increases continuously with compression.
        The regularized thermal switching function is evaluated using
        $\beta=5$ and $\varepsilon_{\theta}=10^{-4}$.
    }
    \label{fig:third_medium_switching_functions}
\end{figure}
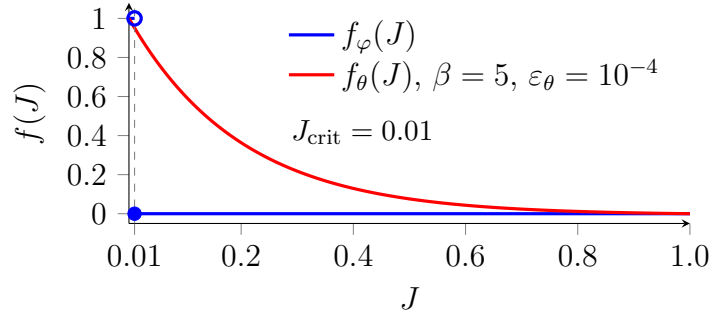

%
%
%
Here, $\beta$ controls the transition width, while $\varepsilon_{\theta}>0$ regularizes the positive-part operator.
For sufficiently small values of $\varepsilon_{\theta}$, the original activation law is recovered within the interval $0\leq J\leq1$, whereas $f_{\theta}(J)$ remains smooth and non-negative for all values of $J$.
Highly compressed regions therefore approach the thermal conductivity of the surrounding solid, while weakly compressed or expanded regions exhibit only limited heat transfer.
Both switching functions are presented in Fig.~\ref{fig:third_medium_switching_functions}.

The switching functions are evaluated in terms of the total volume ratio $J=\det\bF$, rather than the elastic volume ratio $J_{\mathrm e}$.
Hence, transport activation is governed by the actual geometric compression of the third-medium and is not directly affected by the multiplicative separation of thermal expansion.

\subsection{Temperature-dependent electrical conductivity}
\begin{figure}[h]
\centering
\begin{tikzpicture}[x=0.02cm,y=3.5cm,line cap=round,line join=round]

\def\Tzero{293.15}
\def\alphaZero{0.0039}
\def\epsilonSigma{0.05}
\def\ybase{0.2}

\draw[->,thin]
	(250,\ybase)
	--
	(650,\ybase)
	node[right] {$\theta\,[\mathrm{K}]$};

\draw[->,thin]
	(250,\ybase)
	--
	(250,1.12)
	node[above]
	{$\sigma_{\mathrm s}(\theta)/\sigma_{\mathrm s,0}$};

\foreach \x in {300,400,500,600}{
	\draw
		(\x,\ybase+0.015)
		--
		(\x,\ybase-0.015)
		node[below] {\scriptsize $\x$};
}

\foreach \y in {0.4,0.6,0.8,1.0}{
	\draw
		(255,\y)
		--
		(245,\y)
		node[left] {\scriptsize $\y$};
}

\draw[densely dashed,gray]
	(\Tzero,\ybase)
	--
	(\Tzero,1.0);

\draw[densely dashed,gray]
	(250,1.0)
	--
	(\Tzero,1.0);

\node[
	anchor=south west,
	font=\scriptsize
] at (300,1.02)
{$\theta_0=293.15\,\mathrm{K}$};

\draw[
	very thick,
	blue,
	domain=293.15:650,
	samples=150,
	smooth
]
plot
	(
		\x,
		{
		1/max(
			1+\alphaZero*(\x-\Tzero),
			\epsilonSigma
		)
		}
	);

\node[
	anchor=north east,
	font=\scriptsize,
	align=right
] at (640,0.95)
{$\displaystyle
\frac{\sigma_{\mathrm s}(\theta)}
{\sigma_{\mathrm s,0}}
=
\frac{1}
{
\max\!\left[
1+\alpha_0(\theta-\theta_0),
0.05
\right]
}
$};

\node[
	anchor=north east,
	font=\scriptsize,
	align=right
] at (640,0.73)
{$\alpha_0
=
3.9\cdot10^{-3}\,\mathrm{K}^{-1}$};

\end{tikzpicture}
\caption{
Temperature-dependent electrical conductivity of copper according to Eq.~\eqref{eq:sigma_temperature} with $\alpha_0=3.9\times10^{-3}\,\mathrm{K}^{-1}$.
Within the considered temperature range, the solid conductivity decreases monotonically with increasing temperature, while the lower bound of the denominator remains inactive.
}
\label{fig:temperature_dependent_conductivity}
\end{figure}
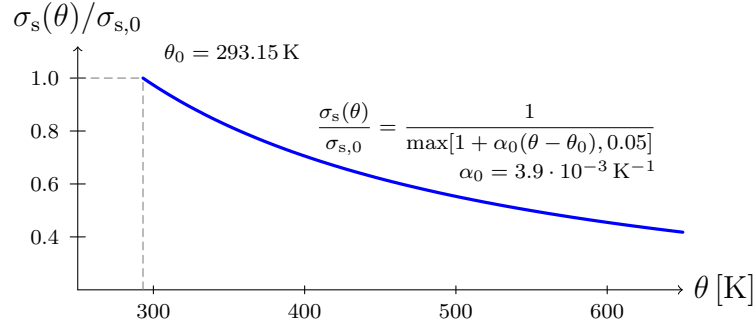

The electrical conductivity of the conducting solid decreases with increasing temperature due to enhanced electron scattering.
Since Joule heating raises the conductor temperature during operation, this effect is included through the temperature-dependent solid conductivity
\begin{equation}
\sigma_{\mathrm s}(\theta)
=
\frac{\sigma_{\mathrm s,0}}
{
\max\!\left[
1+\alpha_0\left(\theta-\theta_0\right),
0.05
\right]
},
\label{eq:sigma_temperature}
\end{equation}
where $\sigma_{\mathrm s,0}$ denotes the electrical conductivity of the solid at the reference temperature $\theta_0$, and $\alpha_0$ is the corresponding temperature coefficient.
A graphical interpretation of these values is given in Fig.~\ref{fig:temperature_dependent_conductivity} for a copper-like material.
The lower bound of the denominator prevents singular or nonphysical conductivity values outside the considered temperature range.
For $\theta\geq\theta_0$, this bound remains inactive and Eq.~\ref{eq:sigma_temperature} reduces to the standard linear resistivity approximation.
Within the third-medium, the conductivity follows as presented in Eq.~\ref{eq:sigma_tm_general}
%
%
such that the temperature dependence is inherited from the conducting solid, while electrical transport is additionally governed by the deformation-dependent switching function $f_{\varphi}(J)$.

\section{The multi-field third-medium finite element implementation}
\vspace{-4mm}

In this section, the finite element implementation of the coupled thermo-electro-mechanical third-medium formulation is presented.
The discretization comprises the displacement field $\bu$, the electric potential $\varphi$, the temperature field $\theta$ and the independently interpolated deformation-gradient-like field $\BTheta$.
For a fixed electrical active set, all primary fields are solved monolithically within the Newton-Raphson scheme.
Changes of the discrete electrical state according to Eq.~\ref{eq:CurrentSwit} are treated by the active-set-type update described in Sec.~\ref{sec:SwitFunc}.
Within an element $\B^{e}$, the primary fields are approximated by
\begin{equation}
\bu
=
\underline{\IN}^{e}_{u}
\,\underline{\bd}^{e}_{u},
\qquad
\varphi
=
\underline{\IN}^{e}_{\varphi}
\,\underline{d}^{e}_{\varphi},
\qquad
\theta
=
\underline{\IN}^{e}_{\theta}
\,\underline{d}^{e}_{\theta},
\qquad\textrm{and}\qquad
\BTheta
=
\underline{\IN}^{e}_{\Theta}
\,\underline{\bd}^{e}_{\Theta},
\label{eq:FE_interpolation_fields}
\end{equation}
where $\underline{\IN}^{e}_{u}$, $\underline{\IN}^{e}_{\varphi}$, $\underline{\IN}^{e}_{\theta}$ and $\underline{\IN}^{e}_{\Theta}$ contain the corresponding shape functions.
The nodal degrees of freedom are collected in $\underline{\bd}^{e}_{u}$, $\underline{d}^{e}_{\varphi}$, $\underline{d}^{e}_{\theta}$ and $\underline{\bd}^{e}_{\Theta}$.
Application of the referential gradient yields
\begin{equation}
\nabla_X\bu
=
\underline{\IB}^{e}_{u}
\,\underline{\bd}^{e}_{u},
\qquad
\nabla_X\varphi
=
\underline{\IB}^{e}_{\varphi}
\,\underline{d}^{e}_{\varphi},
\qquad
\nabla_X\theta
=
\underline{\IB}^{e}_{\theta}
\,\underline{d}^{e}_{\theta},
\qquad\textrm{and}\qquad
\nabla_X\BTheta
=
\underline{\IB}^{e}_{\Theta}
\,\underline{\bd}^{e}_{\Theta}.
\label{eq:FE_interpolation_gradients}
\end{equation}
Consequently, the deformation gradient and the referential electric field are evaluated as
\begin{equation}
\bF
=
\bI
+
\underline{\IB}^{e}_{u}
\,\underline{\bd}^{e}_{u},
\qquad\textrm{and}\qquad
\bE
=
-
\underline{\IB}^{e}_{\varphi}
\,\underline{d}^{e}_{\varphi}.
\label{eq:FE_discrete_F_E}
\end{equation}
For the transient thermal problem, an implicit Euler scheme is employed.
The material time derivative of the temperature is therefore approximated by
$\dot{\theta}_{n+1}
=
\frac{
\theta_{n+1}
-
\theta_{n}
}{
\Delta t
}$,
where $\Delta t$ denotes the time increment.
All element degrees of freedom are collected in
$
\underline{\bd}^{e}
~=~
\begin{bmatrix}
\underline{\bd}^{e}_{u},
\underline{d}^{e}_{\varphi},
\underline{d}^{e}_{\theta},
\underline{\bd}^{e}_{\Theta},
\end{bmatrix}^T.
$
The corresponding virtual fields are interpolated analogously.
For elements belonging to the conducting solid domain, the components associated with $\BTheta$ are omitted.
Insertion of the discrete fields into the weak form gives the element contribution
\begin{equation}
G^{e}
=
G^{e}_{u}
+
G^{e}_{\varphi}
+
G^{e}_{\theta}
+
G^{e}_{\Theta},
\label{eq:FE_element_weak_form}
\end{equation}
where $G^{e}_{\Theta}$ contains the penalty and gradient-regularization contributions of the third-medium.
The element residual vector as well as the system matrix are obtained by differentiation with respect to the virtual and real nodal degrees of freedom as
\begin{equation}
\underline{\IR}^{e}
=
\frac{
\partial
G^{e}
\left(
\underline{\bd}^{e};
\delta\underline{\bd}^{e}
\right)
}{
\partial
\delta\underline{\bd}^{e}
}
\qquad\textrm{and}\qquad
\underline{\IK}^{e}
=
\frac{
\partial
\underline{\IR}^{e}
}{
\partial
\underline{\bd}^{e}
}
\label{eq:FE_element_residual}
\end{equation}
leading to the consistent element system matrix 
\begin{equation}
\underline{\IK}^{e}
=
\begin{bmatrix}
\underline{\IK}^{e}_{uu}
&
\underline{\mathbf 0}
&
\underline{\IK}^{e}_{u\theta}
&
\underline{\IK}^{e}_{u\Theta}
\\[1mm]
\underline{\IK}^{e}_{\varphi u}
&
\underline{\IK}^{e}_{\varphi\varphi}
&
\underline{\IK}^{e}_{\varphi\theta}
&
\underline{\mathbf 0}
\\[1mm]
\underline{\IK}^{e}_{\theta u}
&
\underline{\IK}^{e}_{\theta\varphi}
&
\underline{\IK}^{e}_{\theta\theta}
&
\underline{\mathbf 0}
\\[1mm]
\underline{\IK}^{e}_{\Theta u}
&
\underline{\mathbf 0}
&
\underline{\mathbf 0}
&
\underline{\IK}^{e}_{\Theta\Theta}
\end{bmatrix}.
\label{eq:FE_element_tangent}
\end{equation}
Direct electro-mechanical feedback on the mechanical balance is neglected, such that $\mathbf K_{u\varphi}=\mathbf 0$.
The electric problem nevertheless depends on the deformation through the referential conductivity tensor, resulting in $\mathbf K_{\varphi u}\neq\mathbf 0$.
The block $\underline{\IK}_{\varphi u}$ results from the deformation dependence of the referential electrical
conductivity, while $\underline{\IK}_{\theta u}$ follows from the deformation-dependent heat flux and Joule
source. The blocks $\underline{\IK}_{\varphi\theta}$ and $\underline{\IK}_{\theta\varphi}$ arise from temperature-dependent
conductivity and Joule heating, respectively.
This consistent tangent matrix is generally non-symmetric, resulting from the thermo-mechanical coupling, the temperature-dependent electrical conductivity, the deformation-dependent transport tensors and the Joule-heating contribution.
In particular,
$
\underline{\IK}^{e}_{u\theta}
\neq
\left(
\underline{\IK}^{e}_{\theta u}
\right)^{T}
$
and 
$
\underline{\IK}^{e}_{\varphi\theta}
\neq
\left(
\underline{\IK}^{e}_{\theta\varphi}
\right)^{T}.
$
Consequently, the coupled thermo-electro-mechanical system requires the solution of a non-symmetric linear system within each Newton iteration.
All derivations rely on the strong automatic differentiation engine of AceGen, cf.~\cite{KorWri:2016:aof}.

Linear $T_1$ interpolation is employed for the displacement, electric-potential and temperature fields. 
Within the third-medium domain, each component of the auxiliary field $\BTheta$ is interpolated independently using the same linear shape functions.
However, this approach works for linear, quadratic and higher triangle and quadrilateral finite elements as well, cf.~\cite{VorSchWri:2026:arm}.
All element contributions are evaluated using a six-point Gaussian quadrature rule.

\section{Numerical examples}
\vspace{-4mm}

In this section, the proposed electro-thermo-mechanically coupled TMC formulation is evaluated. 
The focus lies on the performance and coupling behavior of the deformation-dependent electrical and thermal conductivity of the third-medium. 
To this end, different boundary value problems are considered, ranging from academic benchmark tests to more application-oriented examples. 
Unless stated otherwise, copper is used as the conducting solid material. 
The corresponding material parameters are listed in Tab.~\ref{tab:copper_material_parameters_simulation}.

\begin{table}[h]
\centering
\caption{
Material parameters of copper used throughout the numerical examples.
Mechanical properties, density and coefficient of thermal expansion are adopted from \cite{Dav:2001:cac}.
Specific heat capacity is taken from \cite{WhiCol:1984:hco}, thermal conductivity from \cite{TouPowHoKle:1970:tcm}, and electrical conductivity together with its temperature coefficient from \cite{Mat:1979:ero}.
}
\label{tab:copper_material_parameters_simulation}
\begin{tabular}{llll}
\hline
Quantity & Symbol & Value & Unit \\
\hline
Bulk modulus
& $K$
& $1.15 \cdot 10^{5}$
& $\mathrm{N\,mm^{-2}}$ \\

Shear modulus
& $\mu$
& $4.10 \cdot 10^{4}$
& $\mathrm{N\,mm^{-2}}$ \\

Density
& $\rho_0$
& $8.96 \cdot 10^{-6}$
& $\mathrm{kg\,mm^{-3}}$ \\

Specific heat capacity
& $c_0$
& $385$
& $\mathrm{J\,kg^{-1}\,K^{-1}}$ \\

Volumetric heat capacity
& $\rho_0 c_0$
& $3.45 \cdot 10^{-3}$
& $\mathrm{J\,mm^{-3}\,K^{-1}}$ \\

Thermal conductivity
& $\kappa^s$
& $0.401$
& $\mathrm{W\,mm^{-1}\,K^{-1}}$ \\

Electrical conductivity at $\theta_0$
& $\sigma_{s,0}$
& $5.96 \cdot 10^{4}$
& $\mathrm{S\,mm^{-1}}$ \\

Coefficient of thermal expansion
& $\alpha_\theta$
& $16.5 \cdot 10^{-6}$
& $\mathrm{K^{-1}}$ \\

Reference temperature
& $\theta_0$
& $293.15$
& $\mathrm{K}$ \\
\hline
\end{tabular}
\end{table}

All simulations within this contribution rely on the adaptive time-stepping procedure of the package AceGen/AceFEM \cite{KorWri:2016:aof}, if not stated differently.
Each load step is considered as converged for a residual norm $\norm{\underline{\IR}}\leq10^{-10}$.

\subsection{Analytical verification of the electro-thermal material response}
\label{sec:analytical_verification}

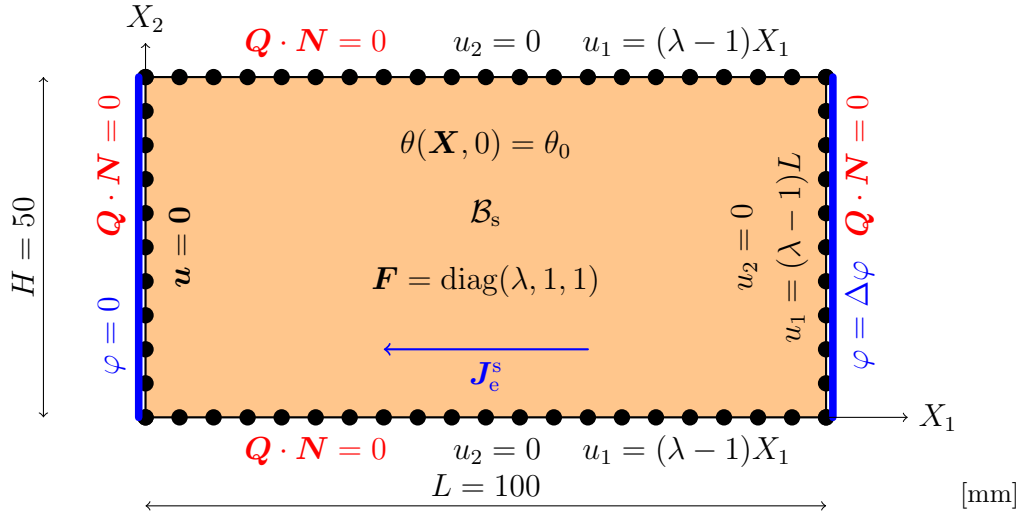
\begin{figure}[h]
\centering

\begin{tikzpicture}[scale=0.09, line cap=round, line join=round]

\def\L{100}
\def\H{50}

\fill[orange!45] (0,0) rectangle (\L,\H);

\draw[thick] (0,0) rectangle (\L,\H);

\node at (50,30) {$\mathcal{B}_{\mathrm{s}}$};

\foreach \y in {0,5,...,50}{
	\fill[black] (0,\y) circle (1.2);
}

\foreach \y in {0,5,...,50}{
	\fill[black] (100,\y) circle (1.2);
}

\foreach \x in {0,5,...,100}{
	\fill[black] (\x,0) circle (1.2);
}

\foreach \x in {0,5,...,100}{
	\fill[black] (\x,50) circle (1.2);
}

\node[black, rotate=90] at (5,25) {$\bu=\bzero$};
\node[black, rotate=90] at (88,25) {$u_2=0$};
\node[black, rotate=90] at (95,25) {$u_1=(\lambda-1)L$};

\node[black] at (70,-5) {$u_2=0$ \quad $u_1=(\lambda-1)X_1$};
\node[black] at (70,55) {$u_2=0$ \quad $u_1=(\lambda-1)X_1$};


\draw[line width=1mm,blue] (-1,0) -- (-1,50);
\draw[line width=1mm,blue] (101,0) -- (101,50);

\node[blue,rotate=90] at (-5,12) {$\varphi=0$};
\node[blue,rotate=90] at (105,15) {$\varphi=\Delta\varphi$};

\node[blue] at (50,6) {$\bJ_{\mathrm e}^{\mathrm s}$};

\draw[->,thick,blue] (65,10) -- (35,10);

\node[red] at (25,55) {$\bQ\cdot\bN=0$};
\node[red] at (25,-5) {$\bQ\cdot\bN=0$};

\node[rotate=90,red] at (-5,37) {$\bQ\cdot\bN=0$};
\node[rotate=90,red] at (105,37) {$\bQ\cdot\bN=0$};

\node at (50,40) {$\theta(\bX,0)=\theta_0$};

\node at (50,20) {$\bF=\mathrm{diag}(\lambda,1,1)$};

\draw[->,thin] (0,0) -- (112,0) node[right] {$X_1$};
\draw[->,thin] (0,0) -- (0,55) node[above] {$X_2$};

\draw[<->] (0,-13) -- (100,-13);
\node at (50,-10) {$L=100$};

\draw[<->] (-15,0) -- (-15,50);
\node[rotate=90] at (-18,25) {$H=50$};

\node[
	anchor=south east,
	font=\footnotesize
] at (130,-15) {[mm]};

\end{tikzpicture}
\caption{
Boundary value problem for the analytical verification of the electro-thermal material response.
A homogeneous deformation $\bF=\mathrm{diag}(\lambda,1,1)$ is prescribed together with an electric potential difference $\Delta\varphi$ across the solid domain $\mathcal{B}_{\mathrm{s}}$.
All thermal boundaries are insulated and the initial temperature is set to $\theta_0$.
}
\label{fig:BVP_validationProblem}
\end{figure}

A homogeneous "1D-like" conducting block, as depicted in Fig.~\ref{fig:BVP_validationProblem}, is considered to verify the electro-thermal constitutive response and its finite element implementation.
Contact effects and deformation-dependent switching functions are excluded from this benchmark.
The mechanical deformation is prescribed by the homogeneous deformation gradient
\begin{equation}
\bF
=
\begin{bmatrix}
\lambda & 0 & 0 \\
0 & 1 & 0 \\
0 & 0 & 1
\end{bmatrix}
\quad\textrm{yielding}\quad
J=\lambda
\quad\textrm{and}\quad
\bC^{-1}
=
\begin{bmatrix}
\lambda^{-2} & 0 & 0 \\
0 & 1 & 0 \\
0 & 0 & 1
\end{bmatrix},
\label{eq:verification_F}
\end{equation}
where $\lambda>0$ denotes the prescribed stretch in the $X_1$-direction.
Thermal expansion is omitted by setting $\alpha_\theta=0$ to isolate the electro-thermal coupling.
A constant electric potential difference is applied over the reference length $L$ such that 
$\varphi(0,t)=0$ and $\varphi(L,t)=\Delta\varphi$ hold.
All thermal boundaries are treated to be insulated as $\bQ\cdot\bN=0$ on $\partial\B$, and the initial temperature is assumed to be spatially homogeneous as $\theta(\bX,0)=\theta_0$.
For the homogeneous deformation and the prescribed electric potential, the exact electric potential field is $\varphi(X_1)=\frac{\Delta\varphi}{L}X_1$.

Accordingly, the referential electric field becomes $\bE=-\nabla_{\!X}\varphi=-\frac{\Delta\varphi}{L}\be_1$.
Using the referential conductivity tensor $\bK_{\mathrm e}^{\mathrm s}=J\sigma_{\mathrm s}(\theta)\bC^{-1},$ from Eq.~\ref{eq:OhmLawReferential} yields
\begin{equation}
\bK_{\mathrm e}^{\mathrm s}
=
\sigma_{\mathrm s}(\theta)
\begin{bmatrix}
\lambda^{-1} & 0 & 0 \\
0 & \lambda & 0 \\
0 & 0 & \lambda
\end{bmatrix}.
\label{eq:verification_conductivity_tensor}
\end{equation}
Hence, the referential current density is spatially constant and follows as $\bJ_{\mathrm e}^{\mathrm s}
=
\bK_{\mathrm e}^{\mathrm s}\cdot\bE
=
-\frac{\sigma_{\mathrm s}(\theta)}{\lambda}
\frac{\Delta\varphi}{L}\be_1$,
with its first component $J_{\mathrm e,1}^{\mathrm s}
=
-\frac{\sigma_{\mathrm s}(\theta)}{\lambda}
\frac{\Delta\varphi}{L}$.
The corresponding Joule heat source is obtained from Eq.~\ref{eq:JouleHeat} as
$Q_{\mathrm J}^{\mathrm s}
=
\bE\cdot\bK_{\mathrm e}^{\mathrm s}\cdot\bE
=
\frac{\sigma_{\mathrm s}(\theta)}{\lambda}
\left(
\frac{\Delta\varphi}{L}
\right)^2$.
Since the temperature remains spatially homogeneous and all boundaries are thermally insulated, the temperature gradient and the heat flux vanish, yielding $\nabla_{\!X}\theta=\bzero$ and $\bQ=\bzero$ and reducing the transient heat balance to
$\rho_0c_0\dot{\theta}
=
\frac{\sigma_{\mathrm s}(\theta)}{\lambda}
\left(
\frac{\Delta\varphi}{L}
\right)^2$.
Within the considered temperature range, the lower bound in Eq.~\ref{eq:sigma_temperature} remains inactive, simplifying the electrical conductivity to 
$\sigma_{\mathrm s}(\theta)
=
\frac{\sigma_{\mathrm s,0}}
{1+\alpha_0(\theta-\theta_0)}$.
Considering both, the transient heat balance and the simplified electrical conductivity leads to 
$\rho_0c_0\dot{\theta}
=
\frac{\sigma_{\mathrm s,0}}
{\lambda\left[1+\alpha_0(\theta-\theta_0)\right]}
\left(
\frac{\Delta\varphi}{L}
\right)^2$.
Introducing the temperature increment $\vartheta(t)=\theta(t)-\theta_0$ and the constant
$
a
=
\frac{\sigma_{\mathrm s,0}}
{\rho_0c_0\lambda}
\left(
\frac{\Delta\varphi}{L}
\right)^2$,
reduces the evolution equation to
$
\dot{\vartheta}
=
\frac{a}{1+\alpha_0\vartheta}$.
The evolution equation
$\dot{\vartheta}=\frac{a}{1+\alpha_0\vartheta}$
can be integrated by separation of variables. 
Using $\dot{\vartheta}=\mathrm{d}\vartheta/\mathrm{d}t$ and multiplying by $1+\alpha_0\vartheta$ gives
$(1+\alpha_0\vartheta)\,\mathrm{d}\vartheta=a\,\mathrm{d}t$.
Integration from the initial state $\vartheta(0)=0$ to the current state $\vartheta(t)$ yields
$\int_0^{\vartheta(t)}(1+\alpha_0\tilde{\vartheta})\,\mathrm{d}\tilde{\vartheta}
=
\int_0^t a\,\mathrm{d}\tau$.
Evaluation of both integrals results in
$\vartheta+\frac{\alpha_0}{2}\vartheta^2=at$.
Solving this quadratic equation for $\vartheta$ and selecting the solution satisfying $\vartheta(0)=0$ gives
$\vartheta(t)=[\sqrt{1+2\alpha_0at}-1]/\alpha_0$.
Consequently, the analytical temperature evolution appears as
\begin{equation}
\theta(t)
=
\theta_0
+
\frac{
\sqrt{1+2\alpha_0at}-1
}{
\alpha_0
}.
\label{eq:verification_temperature_solution}
\end{equation}
A substitution of Eq.~\ref{eq:verification_temperature_solution} into the 1D current density provides the exact transient current density as
\begin{equation}
J_{\mathrm e,1}^{\mathrm s}(t)
=
-\frac{\sigma_{\mathrm s,0}}{\lambda}
\frac{\Delta\varphi}{L}
\frac{1}
{
1+\alpha_0\left[\theta(t)-\theta_0\right]
}.
\label{eq:verification_current_solution}
\end{equation}
This benchmark is discretized using linear $T_1$ finite elements.
Since the exact electric potential and temperature fields are spatially homogeneous or linear, the spatial solution can be represented exactly by the finite element interpolation.
Remaining discrepancies therefore originate from the time discretization.

\begin{Figure}[H]
\unitlength1cm
\begin{picture}(0,5.0)
\put(0.2,5.5){$\lambda=1.0$}
\put(0.2,0.0){\includegraphics[width=7.5cm]{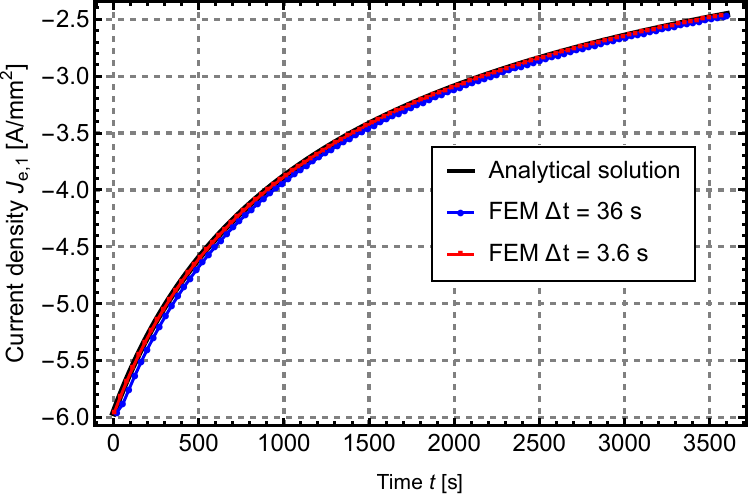}}
\put(8.3,0.0){\includegraphics[width=7.5cm]{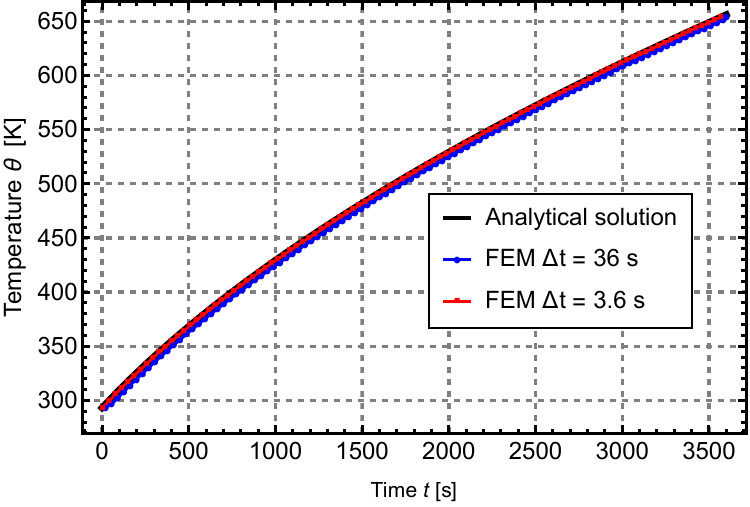}}
\put(0.0,0.1){a)}
\put(8.5,0.1){b)}
\end{picture}
\caption{
Comparison of the analytical solution and the finite element results for $\lambda=1.0$ using $\Delta t=36\,\mathrm{s}$ and $\Delta t=3.6\,\mathrm{s}$:
a) current-density component $J_{\mathrm e,1}$ and
b) temperature $\theta$.
}
\label{fig:AnalyticLambda1}
\end{Figure}

\begin{Figure}[H]
\unitlength1cm
\begin{picture}(0,5.5)
\put(0.2,5.5){$\lambda=1.5$}
\put(0.2,0.0){\includegraphics[width=7.5cm]{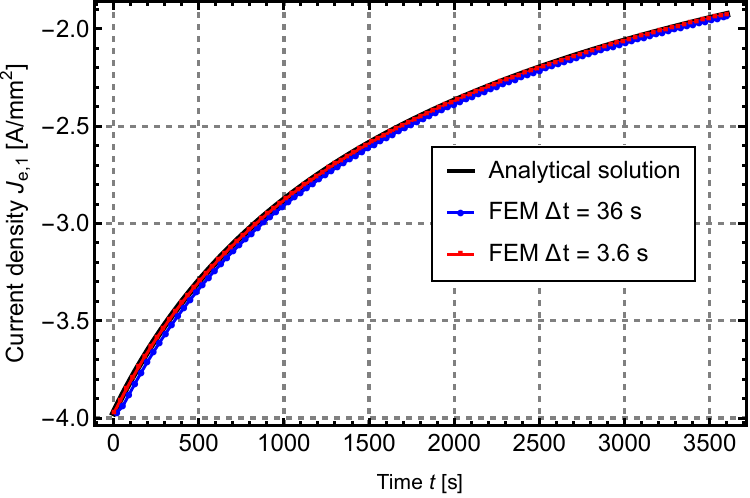}}
\put(8.3,0.0){\includegraphics[width=7.5cm]{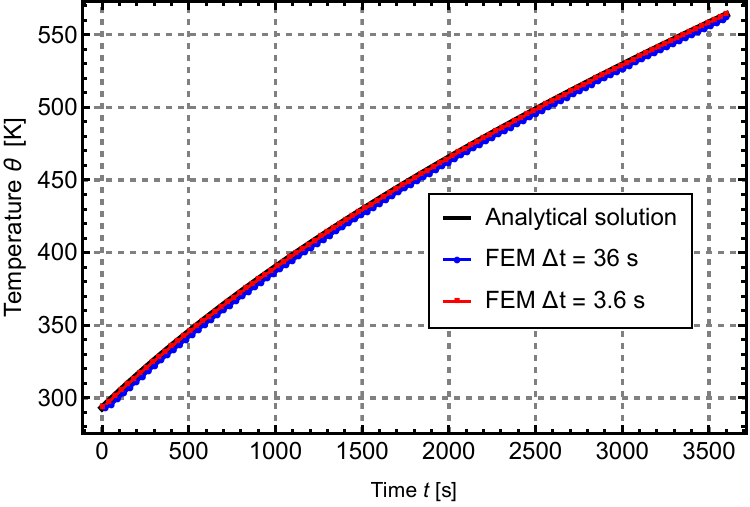}}
\put(0.0,0.1){a)}
\put(8.5,0.1){b)}
\end{picture}
\caption{
Comparison of the analytical solution and the finite element results for $\lambda=1.5$ using $\Delta t=36\,\mathrm{s}$ and $\Delta t=3.6\,\mathrm{s}$:
a) current-density component $J_{\mathrm e,1}$ and
b) temperature $\theta$.
}
\label{fig:AnalyticLambda15}
\end{Figure}

Verification is performed for prescribed stretches $\lambda\in\{1.0,\,1.5\}$, which additionally assesses the deformation dependence of the referential conductivity tensor. 
A time-step refinement study is performed using the time increments $\Delta t\in\{36,3.6\}\,\mathrm{s}$.
The numerical temperature and current density are compared with Eqs.~\ref{eq:verification_temperature_solution} and \ref{eq:verification_current_solution}, respectively.
Agreement with the analytical solution verifies the implemented electro-thermal coupling, including the Piola transformation of the electrical conductivity, Joule heating, temperature-dependent electrical conductivity and transient heat storage.
The numerical results are compared with the analytical solution in Figs.~\ref{fig:AnalyticLambda1} and~\ref{fig:AnalyticLambda15} for $\lambda=1.0$ and $\lambda=1.5$, respectively.
For both deformation states, the finite element solution closely follows the analytical evolution of the current density and temperature over the entire considered time interval.
The remaining deviation is associated with the temporal discretization and decreases consistently when the time-step size is reduced from $\Delta t=36\,\mathrm{s}$ to $\Delta t=3.6\,\mathrm{s}$.
In particular, the results obtained with $\Delta t=3.6\,\mathrm{s}$ are nearly indistinguishable from the analytical reference solutions.
The agreement for both $\lambda=1.0$ and $\lambda=1.5$ verifies the implementation of the deformation-dependent electrical conduction, temperature-dependent conductivity, and resulting Joule heating within the coupled formulation.

\subsection{Contact between two conducting blocks}

\vspace{-5mm}
\begin{figure}[h]
\centering

\begin{minipage}[c]{0.53\textwidth}
\centering
\begin{tikzpicture}[scale=0.05, line cap=round, line join=round]

\def\L{100}
\def\Hs{50}
\def\Htm{50}
\def\Htot{150}

\fill[orange!45] (0,0) rectangle (\L,\Hs);
\fill[yellow!35] (0,\Hs) rectangle (\L,\Hs+\Htm);
\fill[orange!45] (0,\Hs+\Htm) rectangle (\L,\Htot);

\draw[thick] (0,0) rectangle (\L,\Htot);
\draw[thick] (0,\Hs) -- (\L,\Hs);
\draw[thick] (0,\Hs+\Htm) -- (\L,\Hs+\Htm);

\node at (50,25) {$\mathcal{B}_{\mathrm{s}}^{l}$};
\node at (50,75) {$\mathcal{B}_{\mathrm{tm}}$};
\node at (50,125) {$\mathcal{B}_{\mathrm{s}}^{u}$};

\foreach \x in {0,5,...,100}{
	\fill[black] (\x,0) circle (1.2);
	\fill[black] (\x,150) circle (1.2);
}

\node[black, below] at (50,-3) {$\bu=\bzero$};
\node[black, above] at (50,153) {$\bu=\bar \bu$};

\draw[very thick,blue] (0,150) -- (100,150);
\draw[very thick,blue] (0,0) -- (100,0);

\node[blue,right] at (70,159) {$\varphi=\bar\varphi$};
\node[blue,right] at (70,-8) {$\varphi=0$};

\node[rotate=90] at (7,75) {$\bJ_{\mathrm e}\cdot\bN=0$};
\node[rotate=90] at (95,75) {$\bJ_{\mathrm e}\cdot\bN=0$};

\node[rotate=90,red] at (-9,75) {$\bQ\cdot\bN=0$};
\node[rotate=90,red] at (109,75) {$\bQ\cdot\bN=0$};

\draw[->,thin] (0,0) -- (112,0) node[right] {$X$};
\draw[->,thin] (0,0) -- (0,162) node[above] {$Y$};

\draw[<->] (0,-25) -- (100,-25);
\node at (50,-20) {$100$};

\draw[<->] (-30,0) -- (-30,50);
\node[rotate=90] at (-35,25) {$50$};

\draw[<->] (-30,50) -- (-30,100);
\node[rotate=90] at (-35,75) {$50$};

\draw[<->] (-30,100) -- (-30,150);
\node[rotate=90] at (-35,125) {$50$};

\node[
	anchor=south east,
	font=\footnotesize
] at (130,-20) {[mm]};

\end{tikzpicture}

\vspace*{-5mm}
\hspace*{-80mm} a)

\end{minipage}
\hfill
\begin{minipage}[c]{0.45\textwidth}

\vspace*{5mm}
\textbf{Boundary conditions:}\\[2mm]
Displacement: $\bar \bu=[0,-50]^T\,\mathrm{mm}$\\[2mm]
Potential:\hspace*{9.5mm} $\bar\varphi=0.01\,\mathrm{V}$\\[2mm]
Electric sides:\hspace*{4mm} $\bJ_{\mathrm e}\cdot\bN=0$\\[2mm]
Thermal sides:\hspace*{3mm} $\bQ\cdot\bN=0$

\vspace*{3mm}
\textbf{Thermal cases:}\\[2mm]
Case 1: $\bQ\cdot\bN=0$ on $\partial\mathcal{B}$\\[1mm]
Case 2: $\theta(Y=0\wedge Y=150)=\theta_0$ 

\vspace*{2mm}
\textbf{Displacement parameter:}

\begin{tikzpicture}[x=1.6cm,y=2.5cm]
\draw[->] (0,0) -- (3.2,0) node[right] {$t$};
\draw[->] (0,0) -- (0,1.1) node[above] {$\lambda$};
\foreach \x/\lab in {0.5/0.5,1/1,3/3600}{
	\draw (\x,0.02) -- (\x,-0.02) node[below] {\scriptsize $\lab$};
}
\foreach \y in {0.2,0.4,0.6,0.8,1}{
	\draw (0.02,\y) -- (-0.02,\y) node[left] {\scriptsize $\y$};
}
\draw[very thick,blue] (0,0) -- (1,1) -- (3,1);
\end{tikzpicture}

\vspace*{-10mm}
b)
\vspace*{8mm}
\end{minipage}

\caption{
Boundary value problem for an idealized high-current switch contact.
Two conducting solid blocks are separated by a third-medium layer and compressed by a prescribed vertical displacement.
A potential difference drives electric current through the system once a sufficiently compressed conducting path forms.
a) shows the geometry and the mechanical, electrical and thermal boundary conditions, while
b) summarizes the loading history and the two thermal cases considered below.
}
\label{fig:bvp_two_blocks_with_conductivity}
\end{figure}
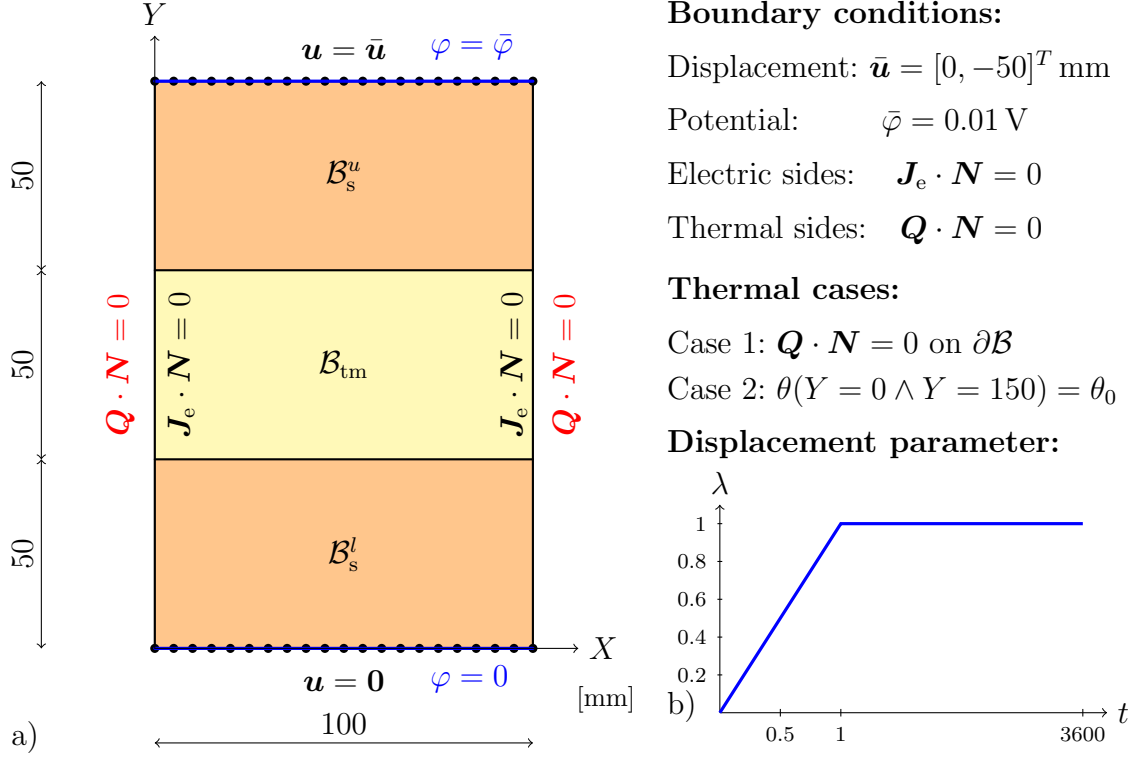

An idealized high-current switch contact is considered.
The conducting solids $\B_{\mathrm{s}}^{u}$ and $\B_{\mathrm{s}}^{l}$ are separated by the third-medium domain $\B_{\mathrm{tm}}$, see Fig.~\ref{fig:bvp_two_blocks_with_conductivity}.
Each domain has a width of $100\,\mathrm{mm}$, while both conducting blocks and the third-medium layer have a height of $50\,\mathrm{mm}$.
The lower boundary is fixed by prescribing $\bu=\bzero$, whereas $\bar u_1=0\,\mathrm{mm}$ and $\bar u_2=-50\,\mathrm{mm}$ are imposed along the upper boundary.
Consequently, the upper conductor approaches the lower one and progressively compresses the intermediate third-medium.
An electric potential difference is applied between the upper and lower conductor surfaces.
The lower boundary is grounded by prescribing $\varphi=0$, whereas the upper boundary is set to $\bar\varphi=0.01\,\mathrm{V}$.
Electrical insulation is imposed along the lateral boundaries through $\bJ_{\mathrm e}\cdot\bN=0$.
Initially, the third-medium is electrically insulating and prevents current transfer between the conductors.
Progressive compression reduces the local volume ratio of the third-medium until electrically conducting regions are activated according to Eq.~\ref{eq:CurrentSwit} and a continuous conducting path is established.
Joule heat generated by the resulting current is redistributed by Fourier heat conduction.
Two thermal boundary conditions are investigated:\\[1mm]
 \textbf{Case 1:}
All external boundaries are thermally insulated by prescribing
$\bQ\cdot\bN=0$.
Accordingly, no heat is exchanged with the surroundings, and the temperature evolution is governed solely by the internal Joule heating and the prescribed initial temperature $\theta_0$.
\textbf{Case 2:}
The upper and lower conductor ends are maintained at $\theta=\theta_0$ and represent ideal heat reservoirs.
The lateral boundaries remain thermally insulated according to $\bQ\cdot\bN=0$.

\begin{Figure}[h]
\begin{picture}(0,13.5)
\unitlength1cm
\put(2.5,8.5){\includegraphics[width=11cm]{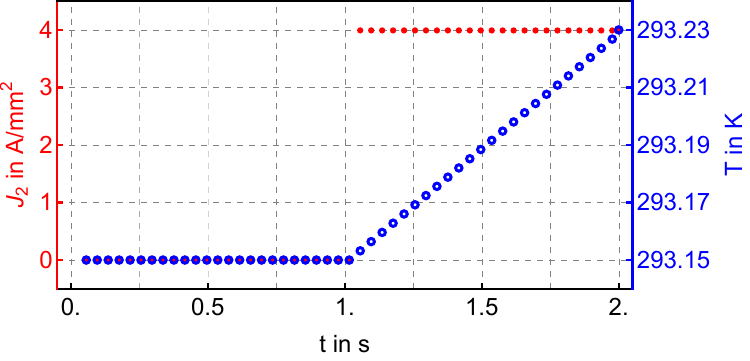}}

\put(0.0,1.1){\includegraphics[width=0.31\textwidth]{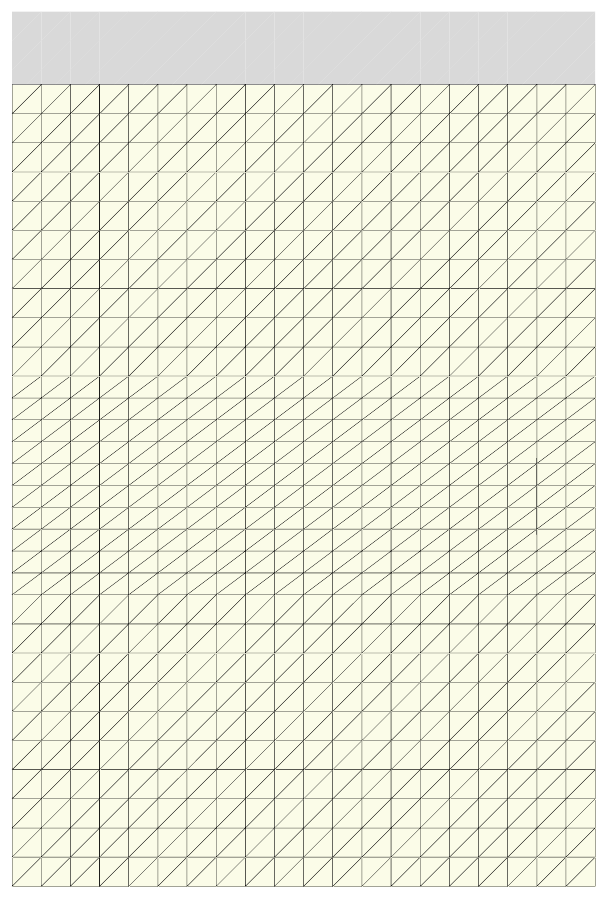}}
\put(5.3,1.1){\includegraphics[width=0.31\textwidth]{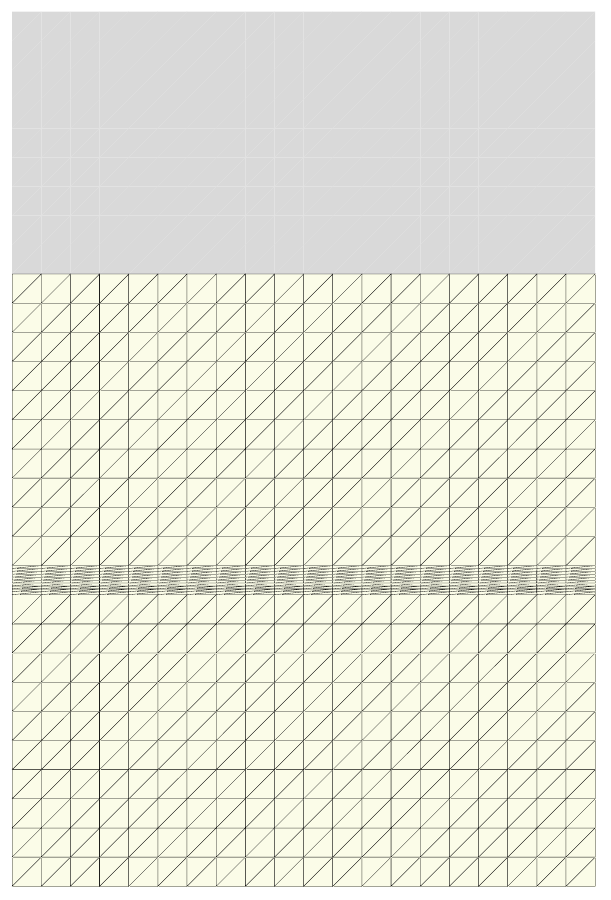}}
\put(10.6,1.1){\includegraphics[width=0.31\textwidth]{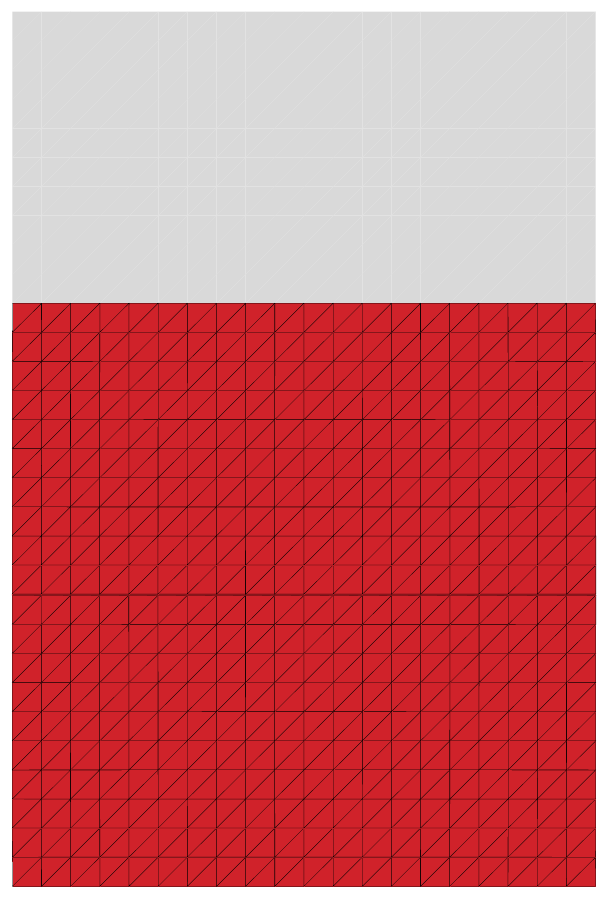}}

\put(2.0,9.0){a)}

\put(0.2,0.4){b)}
\put(5.5,0.4){c)}
\put(10.8,0.4){d)}

\put(1.9,0.4){$\lambda=0.25$}
\put(7.1,0.4){$\lambda=0.9$}
\put(12.5,0.4){$\lambda=1$}

\put(4.0,10.5){\textcolor{red}{b)}}
\put(7.3,10.5){\textcolor{red}{c)}}
\put(7.7,10.5){\textcolor{red}{d)}}

\end{picture}
\vspace*{-5mm}
\caption{
Initial closing process for Case~1.
a) shows the time histories of the current-density component $J_{\mathrm e,2}$ and the temperature $\theta$ over $t\in[0,2]\,\mathrm{s}$, while
b)-d) present the corresponding contact configurations at $\lambda=0.25$, $\lambda=0.9$ and $\lambda=1$, respectively.
Conducting contact is established at $\lambda=1$.
}
\label{fig:TempCond}
\end{Figure}

Geometry, mechanical loading and electrical boundary conditions remain identical in both cases.
The computational domain is discretized using 1200 linear $T_1$ finite elements.
Third-medium parameters are chosen as $p_\Theta=10^4$, $\alpha_r=10^2$ and $\gamma=10^{-6}$ and remain unchanged throughout this benchmark.
Fig.~\ref{fig:TempCond} illustrates the initial closing process for Case~1.
Up to $t=1\,\mathrm{s}$, the conductors remain electrically disconnected and no current is conducted.
Consequently, no Joule heat is generated and the temperature remains at its initial value.
Once a continuous conducting path is activated by the switching criterion, the electric circuit closes and the current density exhibits the corresponding abrupt onset imposed by the discontinuous electrical activation law.
In contrast, the temperature remains continuous because its evolution is governed by the transient heat balance.
Since no heat is removed through the external boundaries in Case 1, the generated Joule heat remains inside the computational domain and produces a monotonic temperature increase.
At $\lambda=0.25$, the conductors are separated and no conducting path exists.
Only a narrow gap remains at $\lambda=0.9$, while electrical conduction is activated at $\lambda=1$.
The results reproduce the expected transition from an open to a closed electrical contact and illustrate the different temporal characteristics of current transfer and heat accumulation.

\begin{Figure}[H]
\unitlength1cm
\begin{picture}(0,4.0)
\put(0.2,0.0){\includegraphics[width=8cm]{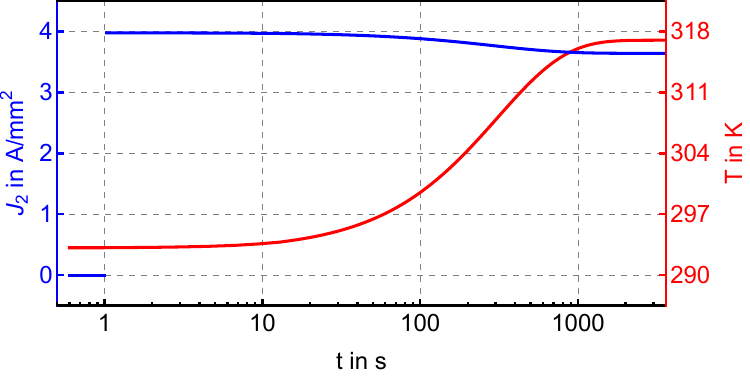}}
\put(8.5,0.0){\includegraphics[width=8cm]{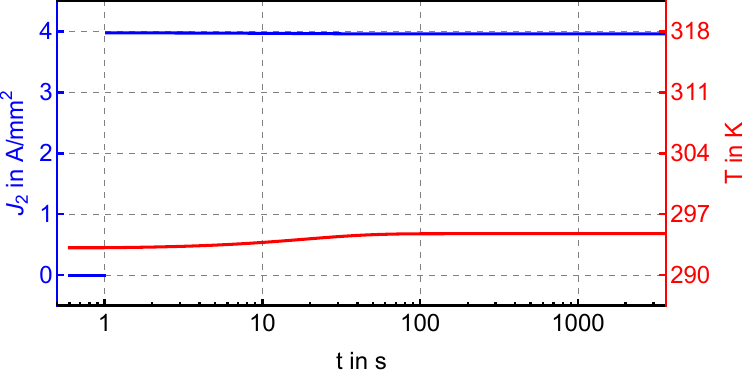}}
\put(0.0,0.1){a)}
\put(8.5,0.1){b)}
\end{picture}
\caption{
Time histories of the current-density component $J_{\mathrm e,2}$ and the temperature $\theta$ on a logarithmic time scale for
a) Case~1 with fully insulated thermal boundaries and
b) Case~2 with isothermal upper and lower boundaries.
Evaluated at point P(50,50).
}
\label{fig:TempCurrentLog}
\end{Figure}

Following the initial closing process, the third-medium remains in the compressed configuration reached at $t=1\,\mathrm{s}$.
The subsequent electro-thermo-mechanical response is evaluated over $I_{\textrm{time}}=[0,3600\,\mathrm{s}]$, corresponding to one hour of operation.
Comparison of Cases~1 and 2 isolates the influence of heat removal on the temperature-dependent electrical conductivity and current density.
In Fig.~\ref{fig:TempCurrentLog} the electro-thermal response of both cases is compared. 
Each data set is evaluated at point P(50,50).
For Case~1, the absence of heat removal through the external boundaries causes a pronounced temperature increase.
The corresponding reduction in electrical conductivity decreases the current density despite the constant applied potential difference.
Hence, the formulation captures the feedback between Joule heating, temperature and electrical conductivity.
For Case~2, heat is removed through the isothermal upper and lower boundaries.
The temperature increase remains moderate and approaches a bounded state.
Consequently, changes in the electrical conductivity are small and the current density remains nearly constant.
The comparison highlights the decisive influence of thermal boundary conditions on the predicted electro-thermal response.

\subsection{Contact between two conducting blocks and non-symmetric loading}
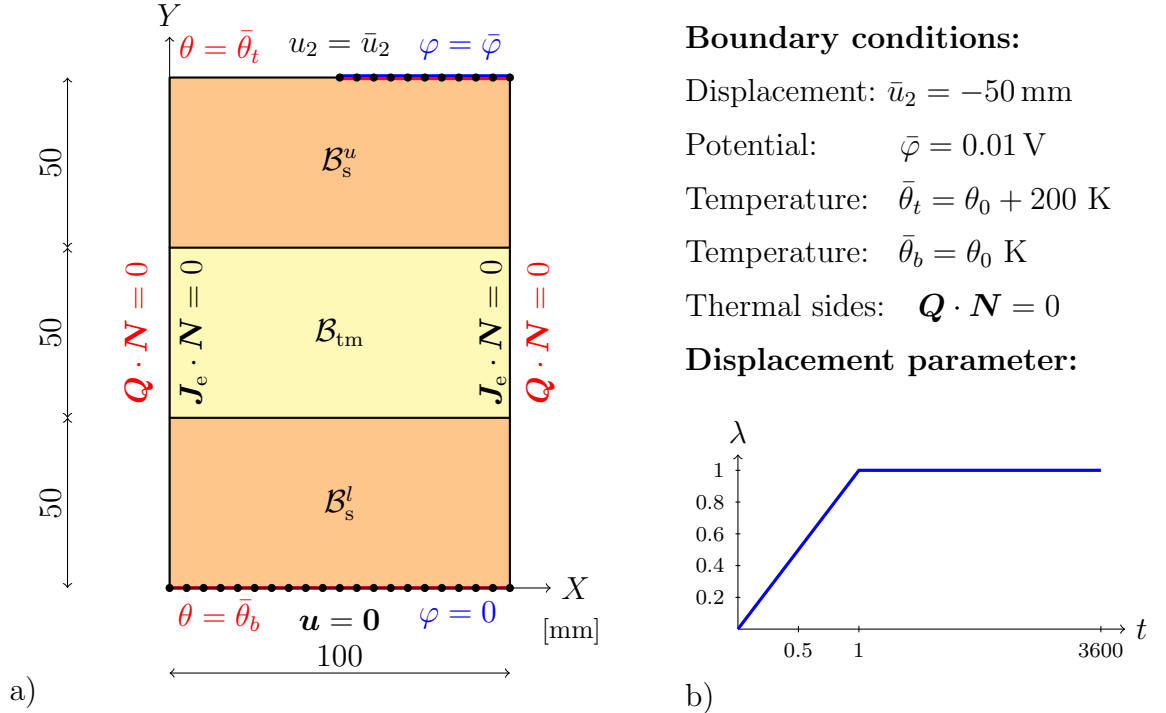
\begin{figure}[h]
\centering

\begin{minipage}[c]{0.58\textwidth}
\centering
\begin{tikzpicture}[scale=0.045, line cap=round, line join=round]

\def\L{100}
\def\Hs{50}
\def\Htm{50}
\def\Htot{150}

\fill[orange!45] (0,0) rectangle (\L,\Hs);
\fill[yellow!35] (0,\Hs) rectangle (\L,\Hs+\Htm);
\fill[orange!45] (0,\Hs+\Htm) rectangle (\L,\Htot);

\draw[thick] (0,0) rectangle (\L,\Htot);
\draw[thick] (0,\Hs) -- (\L,\Hs);
\draw[thick] (0,\Hs+\Htm) -- (\L,\Hs+\Htm);

\node at (50,25) {$\mathcal{B}_{\mathrm{s}}^{l}$};
\node at (50,75) {$\mathcal{B}_{\mathrm{tm}}$};
\node at (50,125) {$\mathcal{B}_{\mathrm{s}}^{u}$};

\draw[very thick, blue] (50,150.5) -- (100,150.5);
\draw[very thick, blue] (0,0) -- (100,0);

\node[blue, right] at (70,159) {$\varphi=\bar\varphi$};
\node[blue, right] at (70,-8) {$\varphi=0$};

\draw[very thick, red] (50,149.7) -- (100,149.7);
\draw[very thick, red] (0,0) -- (100,0);

\node[red, left] at (30,159) {$\theta=\bar \theta_t$};
\node[red, left] at (30,-8) {$\theta=\bar \theta_b$};

\node[rotate=90,red] at (-9,75) {$\bQ\cdot\bN=0$};
\node[rotate=90,red] at (109,75) {$\bQ\cdot\bN=0$};

\node[rotate=90] at (7,75) {$\bJ_{\mathrm e}\cdot\bN=0$};
\node[rotate=90] at (95,75) {$\bJ_{\mathrm e}\cdot\bN=0$};

\foreach \x in {0,5,...,100}{
	\fill[black] (\x,0) circle (1.2);
}

\foreach \x in {50,55,...,100}{
	\fill[black] (\x,150) circle (1.2);
}

\node[black, below] at (50,-3) {$\bu=\bzero$};
\node[black, above] at (50,153) {$u_2=\bar u_2$};

\draw[->, thin] (0,0) -- (112,0) node[right] {$X$};
\draw[->, thin] (0,0) -- (0,162) node[above] {$Y$};

\draw[<->] (0,-25) -- (100,-25);
\node at (50,-20) {$100$};

\draw[<->] (-30,0) -- (-30,50);
\node[rotate=90] at (-35,25) {$50$};

\draw[<->] (-30,50) -- (-30,100);
\node[rotate=90] at (-35,75) {$50$};

\draw[<->] (-30,100) -- (-30,150);
\node[rotate=90] at (-35,125) {$50$};

\node[
    anchor=south east,
    font=\footnotesize
] at (130,-20) {[mm]};

\end{tikzpicture}

\vspace*{-2mm}
\hspace*{-80mm} a)

\end{minipage}
\hfill
\begin{minipage}[c]{0.40\textwidth}

\vspace*{5mm}
\textbf{Boundary conditions:}\\[2mm]
Displacement: $\bar u_2=-50\,\mathrm{mm}$\\[2mm]
Potential:\hspace*{9.5mm} $\bar\varphi=0.01\,\mathrm{V}$\\[2mm]
Temperature:\hspace*{2.5mm} $\bar \theta_t=\theta_0+200$ K\\[2mm]
Temperature:\hspace*{2.5mm} $\bar \theta_b=\theta_0$ K\\[2mm]
Thermal sides:\hspace*{3mm} $\bQ\cdot\bN=0$

\vspace*{2mm}
\textbf{Displacement parameter:}\\

\begin{tikzpicture}[x=1.6cm,y=2.1cm]
\draw[->] (0,0) -- (3.2,0) node[right] {$t$};
\draw[->] (0,0) -- (0,1.1) node[above] {$\lambda$};
\foreach \x/\lab in {0.5/0.5,1/1,3/3600}{
    \draw (\x,0.02) -- (\x,-0.02) node[below] {\scriptsize $\lab$};
}
\foreach \y in {0.2,0.4,0.6,0.8,1}{
    \draw (0.02,\y) -- (-0.02,\y) node[left] {\scriptsize $\y$};
}
\draw[very thick,blue] (0,0) -- (1,1) -- (3,1);
\end{tikzpicture}
\vspace*{-12mm}b)
\vspace*{12mm}
\end{minipage}
\vspace*{-4mm}
\caption{
Boundary value problem for localized electro-thermo-mechanical contact. 
Displacement, electric potential and temperature are prescribed on the right half of the upper boundary, while the lower boundary conditions remain spatially uniform.
Both lateral boundaries are electrically and thermally insulated.}
\label{fig:bvp_two_blocks_with_conductivity_skew}
\end{figure}

In a second example, the previously considered boundary value problem is extended by introducing a localized loading configuration, see Fig.~\ref{fig:bvp_two_blocks_with_conductivity_skew}. 
Instead of prescribing the mechanical displacement, electric potential and temperature along the complete upper boundary, all three boundary conditions are applied only to the right half of the upper conductor.

The lower boundary conditions remain unchanged. 
Consequently, contact is established locally and the electric current enters the conducting domain only through the loaded part of the interface. 
To introduce an additional thermal driving force, the prescribed temperature at the upper boundary is increased to $\theta_0+200\,\mathrm{K}$, whereas the lower boundary is maintained at the reference temperature $\theta_0$. 
The remaining part of the upper boundary is mechanically traction-free,
electrically and thermally insulated.

\begin{Figure}[H]
\begin{picture}(0,5.5)
\unitlength1cm
\put( 1.1,0.5){\includegraphics[width=4.cm]{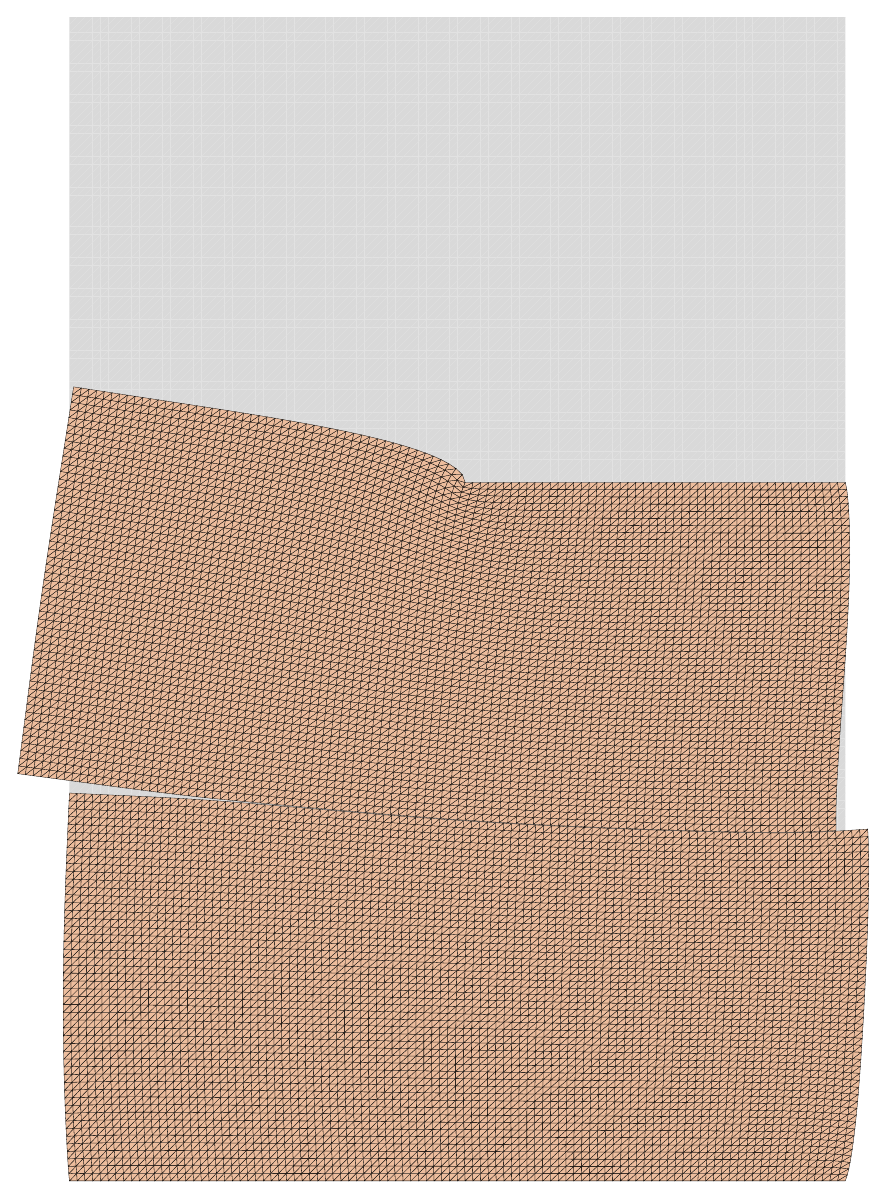}}
\put( 5.8,0.5){\includegraphics[width=5cm]{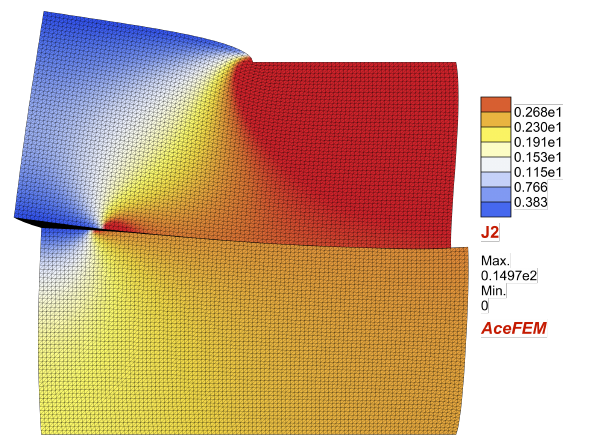}}
\put(10.5,0.5){\includegraphics[width=5cm]{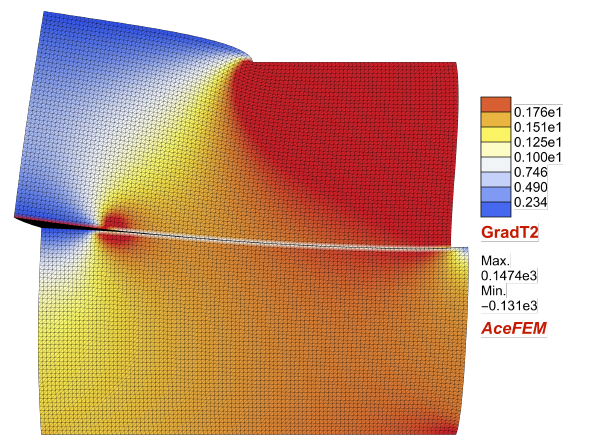}}
\put( 1.0,0.5){a)}
\put( 5.7,0.5){b)}
\put(10.4,0.5){c)}
\end{picture}
\vspace*{-5mm}
\caption{
Results for the localized thermo-electro-mechanical contact problem. 
a) shows the deformed configuration after contact formation, 
b) depicts the resulting current density $J_{e,2}$, illustrating the localized current path through the contact region,
c) presents the corresponding temperature field.
All quantities are shown at $t=3600\,\mathrm{s}$ after the system has approached a quasi-steady thermal state.
}
\label{fig:TempCondSkew}
\end{Figure}

The lateral boundaries remain electrically and thermally insulated. 
Compared with the previous benchmark, this setup produces strongly non-uniform mechanical, electrical and thermal fields and therefore represents a considerably more demanding demonstration problem for the proposed formulation.
The computational domain is discretized by approximately $3\times10^4$ linear triangular elements. 
The transient response is evaluated over $3600\,\mathrm{s}$ using 192 load increments, which require a total of 984 Newton iterations.
Fig.~\ref{fig:TempCondSkew} demonstrates the fully coupled response for the localized loading configuration. 
The imposed displacement leads to contact formation pronounced on the right side of the interface, resulting in a strongly non-uniform current distribution. 
Accordingly, the highest current densities are concentrated in the contact region and decay continuously towards the unloaded part of the interface. 
Joule heating follows the same localization and produces elevated temperatures in the vicinity of the conducting contact. 
The prescribed temperature difference between the upper and lower boundaries generates an additional heat flux through the structure, which lets the temperature field decay from the top via the contact region to the bottom. 
The corresponding temperature gradient is largest close to the localized heat source and gradually decreases towards the remaining parts of the domain. 
Despite the strongly coupled and highly localized fields, all solution variables remain smooth and physically consistent, demonstrating the workability of the proposed electro-thermo-mechanical third-medium formulation.
In Fig.~\ref{fig:TempCondSkewPeriod} the time evolution of the temperature field after the conducting contact has been established is illustrated. 
Immediately after contact formation, elevated temperatures are confined to the prescribed hot boundary and to the neighboring contact region, where Joule heating is concentrated. 
As time progresses, thermal diffusion transports the generated heat into the remaining parts of both conducting blocks, resulting in an increasingly smooth temperature distribution. 
The strongest temperature gradients are observed during the early stages of the simulation and gradually decrease as the heat penetrates the entire domain. 
After approximately $100\,\mathrm{s}$, the transient behavior becomes noticeably slower, indicating that thermal diffusion dominates over the initially localized heat generation. 
At $t=3600\,\mathrm{s}$, the solution approaches a quasi-steady temperature field with a smooth spatial distribution between the prescribed hot upper boundary and the colder lower boundary. 
The results demonstrate the consistent interaction between localized Joule heating, prescribed thermal boundary conditions, and transient heat conduction within the coupled electro-thermo-mechanical formulation.

%
%
%
%
%
%
%
%
%
%

\begin{Figure}[H]
\begin{picture}(0,12.0)
\unitlength1cm

\put( 0.0,8.0){\includegraphics[width=5.50cm]{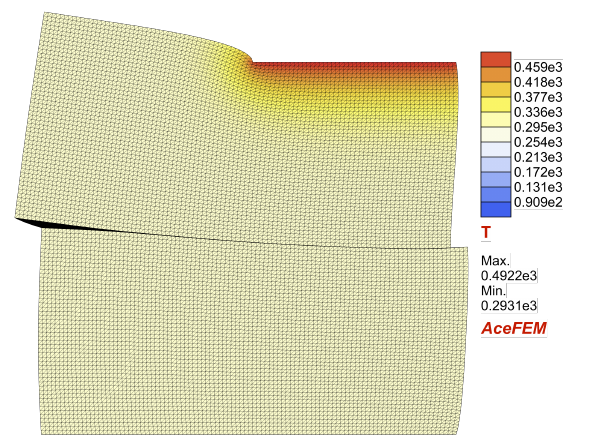}}
\put( 5.7,8.0){\includegraphics[width=5.50cm]{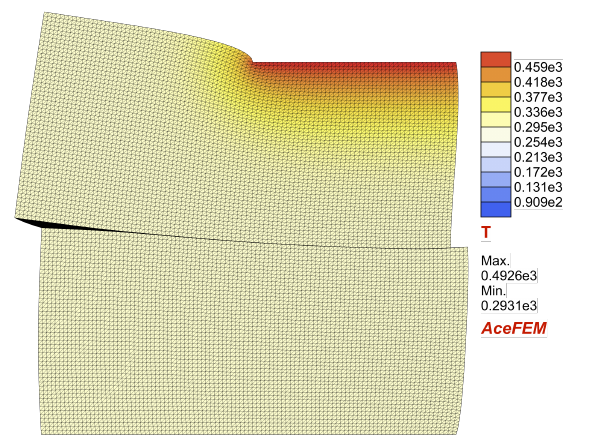}}
\put(11.0,8.0){\includegraphics[width=5.50cm]{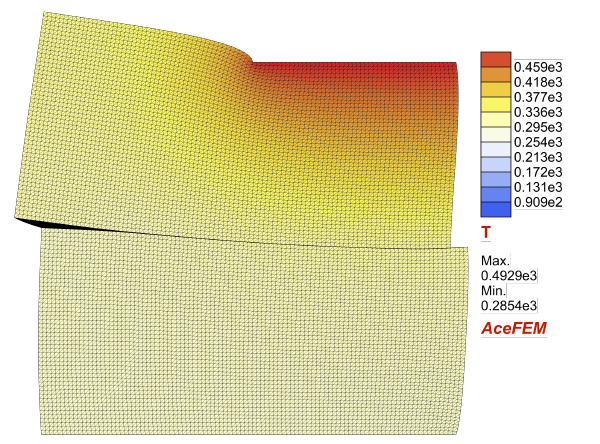}}

\put( 0.0,3.5){\includegraphics[width=5.50cm]{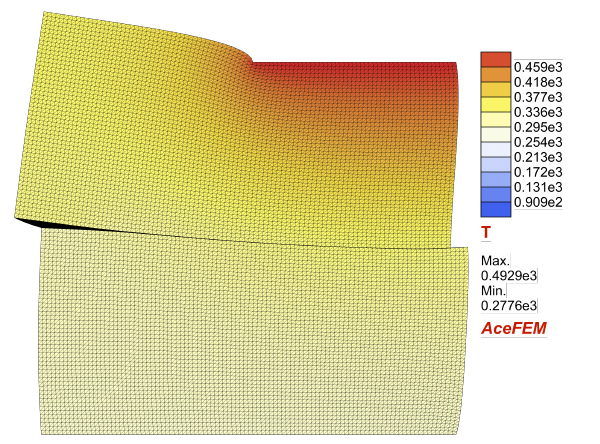}}
\put( 5.7,3.5){\includegraphics[width=5.50cm]{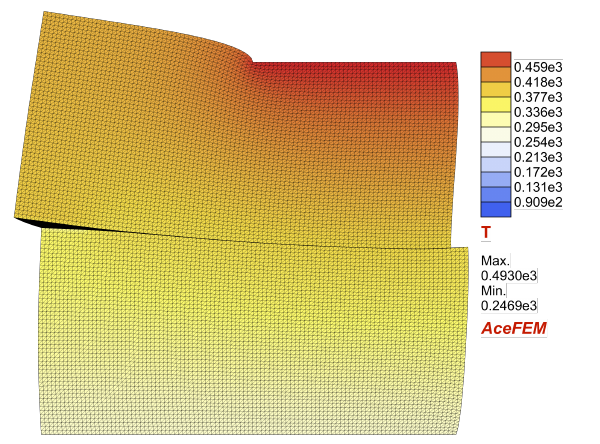}}
\put(11.0,3.5){\includegraphics[width=5.50cm]{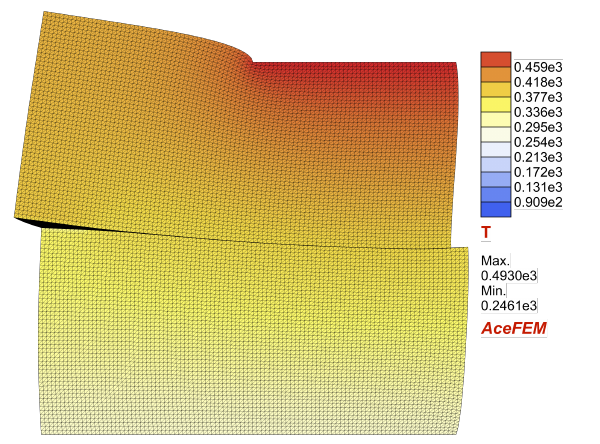}}

\put( 0.0,7.9){a)}
\put( 5.6,7.9){b)}
\put(10.8,7.9){c)}

\put( 0.0,3.2){d)}
\put( 5.6,3.2){e)}
\put(10.8,3.2){f)}

\end{picture}
\vspace*{-35mm}
\caption{
Transient evolution of the temperature field after contact closure at
a)  1.1 s,
b)  2.0 s,
c) 10.0 s,
d) 20.0 s,
e) 100.0 s,
f) 3600 s.
}
\label{fig:TempCondSkewPeriod}
\end{Figure}

\subsection{Contact between blocks with microstructural surface roughness}

\begin{figure}[H]
\centering

\begin{minipage}[c]{0.6\textwidth}
\centering

\begin{tikzpicture}[scale=0.045, line cap=round, line join=round]

\def\L{100}
\def\H{125}
\def\A{7}

%
%
%


\fill[orange!45]
    (0,0)
    --
    (\L,0)
    --
    plot[
        domain=\L:0,
        samples=300,
        variable=\x
    ]
    ({\x},{44+\A*cos(60+10.8*\x)})
    --
    cycle;

\fill[orange!45]
    (0,\H)
    --
    (\L,\H)
    --
    plot[
        domain=\L:0,
        samples=300,
        variable=\x
    ]
    ({\x},{79+\A*sin(10.8*\x)})
    --
    cycle;

\fill[yellow!35]
    plot[
        domain=0:\L,
        samples=300,
        variable=\x
    ]
    ({\x},{44+\A*cos(60+10.8*\x)})
    --
    plot[
        domain=\L:0,
        samples=300,
        variable=\x
    ]
    ({\x},{79+\A*sin(10.8*\x)})
    --
    cycle;

\draw[thick]
    (0,0)
    rectangle
    (\L,\H);

\draw[thick]
    plot[
        domain=0:\L,
        samples=300,
        variable=\x
    ]
    ({\x},{44+\A*cos(60+10.8*\x)});

\draw[thick]
    plot[
        domain=0:\L,
        samples=300,
        variable=\x
    ]
    ({\x},{79+\A*sin(10.8*\x)});

\node at (50,18)
    {$\mathcal{B}_{\mathrm{s}}$};

\node at (50,62)
    {$\mathcal{B}_{\mathrm{tm}}$};

\node at (50,105)
    {$\mathcal{B}_{\mathrm{s}}$};


\draw[very thick, blue]
    (0,1.0)
    --
    (100,1.0);

\node[blue, right]
    at (68,-8)
    {$\varphi=0$};

\draw[very thick, blue]
    (00,126.)
    --
    (100,126.);

\node[blue, right]
    at (68,134)
    {$\varphi=\bar\varphi$};

\node[rotate=90]
    at (-9,90)
    {$\bJ_{\mathrm e}\cdot\bN=0$};

\node[rotate=90]
    at (109,90)
    {$\bJ_{\mathrm e}\cdot\bN=0$};


\draw[very thick, red]
    (0,-0.5)
    --
    (100,-0.5);

\node[red, left]
    at (32,-8)
    {$\theta=\bar \theta_b$};

\draw[very thick, red]
    (0,124.)
    --
    (100,124.);

\node[red, left]
    at (32,134)
    {$\theta=\bar \theta_t$};

\node[rotate=90, red]
    at (-9,35)
    {$\bQ\cdot\bN=0$};

\node[rotate=90, red]
    at (109,35)
    {$\bQ\cdot\bN=0$};


\foreach \x in {0,5,...,100}{
    \fill[black]
        (\x,0)
        circle (1.5);
}

\foreach \x in {0,5,...,100}{
    \fill[black]
        (\x,125)
        circle (1.5);
}

\foreach \y in {3,6,...,122}{
    \fill[cyan]
        (0,\y)
        circle (1.2);

    \fill[cyan]
        (100,\y)
        circle (1.2);
}

\node[black, below]
    at (50,-3)
    {$\bu=\bzero$};

\node[black, above]
    at (50,128)
    {$\bu=\bar\bu$};

\node[cyan, rotate=90]
    at (6,20)
    {$u_1=0$};

\node[cyan, rotate=90]
    at (94,20)
    {$u_1=0$};

\draw[->, thin]
    (0,0)
    --
    (112,0)
    node[right] {$X$};

\draw[->, thin]
    (0,0)
    --
    (0,137)
    node[above] {$Y$};


\draw[<->]
    (0,-25)
    --
    (100,-25);

\node at (50,-20)
    {$2$};

\draw[<->]
    (-24,0)
    --
    (-24,125);

\node[rotate=90]
    at (-29,62.5)
    {$2.5$};

\node[anchor=east]
    at (-2,125)
    {$2$};

\node[
    anchor=south east,
    font=\footnotesize
]
at (124,-18)
{$[\mathrm{mm}]$};

\end{tikzpicture}

\end{minipage}
\hspace*{-5mm}
\begin{minipage}[c]{0.40\textwidth}

\vspace*{5mm}
\textbf{Boundary conditions:}\\[2mm]
Displacement: $\bar \bu=[0,-2]^T\,\mathrm{mm}$\\[2mm]
Potential:\hspace*{9.5mm} $\bar\varphi=0.01\,\mathrm{V}$\\[2mm]
Temperature:\hspace*{2.5mm} $\bar \theta_t=\theta_0+200$ K\\[2mm]
Temperature:\hspace*{2.5mm} $\bar \theta_b=\theta_0$ K\\[2mm]
\vspace*{4mm}

\end{minipage}

\caption{
Boundary value problem for the compression of a solid layer with rough interfaces.
}
\label{fig:bvp_rough_thermo_electro_mechanical}

\end{figure}

Contact between nominally rough surfaces is governed by the interaction of individual asperities rather than by the apparent contact area.
The resulting constriction of electric current strongly affects the effective electrical contact resistance.
A geometrically idealized rough-surface configuration, inspired by the rough-interface example of Wriggers~\cite{Wri:2026:atm}, is therefore considered to assess whether the proposed third-medium formulation captures progressive gap closure and the associated localization of electrical transport.
The example is intended as a qualitative demonstration of the proposed electrical transport formulation rather than as a quantitative model of electrical contact resistance at real rough metallic interfaces.
Two solid domains are separated by a deformable third-medium layer that represents the initial gap between the rough surfaces.
The lower and upper interfaces are prescribed by
$
y_{\mathrm{l}}(x)
=
0.38
+
0.14
\cos\left(
\frac{\pi}{3}
+
\frac{6\pi x}{L}
\right)
$
and
$
y_{\mathrm{u}}(x)
=
1.08
+
0.14
\sin\left(
\frac{6\pi x}{L}
\right),
$
respectively, where $L=2$ mm denotes the specimen width, and both surfaces contain three asperity periods across the specimen width.
The lower boundary is fully constrained, whereas
$\bar{\bu}=[0,-2]^{\mathsf T}\,\mathrm{mm}$
is prescribed along the upper boundary.
In addition, $u_1=0$ is enforced on both lateral sides.
Electric potentials $\varphi=0$ and $\varphi=\bar{\varphi}=0.01\,\mathrm{V}$ are applied at the lower and upper boundaries, respectively.
Temperatures $\theta=\theta_0$ and $\theta=\theta_0+200\,\mathrm{K}$ are prescribed at the lower and upper boundaries, while both lateral sides are electrically and thermally insulated.
The complete domain is discretized by $17\,664$ linear triangular $T_1$ elements.
The coupled problem is solved over $t\in[0,1]\,\mathrm{s}$ using adaptive time-step control.
In total, $247$ accepted load steps and $1\,957$ Newton iterations are required.

\begin{Figure}[H]
\begin{picture}(0,6.5)
\unitlength1cm

\put( 0.7,0.8){\includegraphics[width=4.7cm]{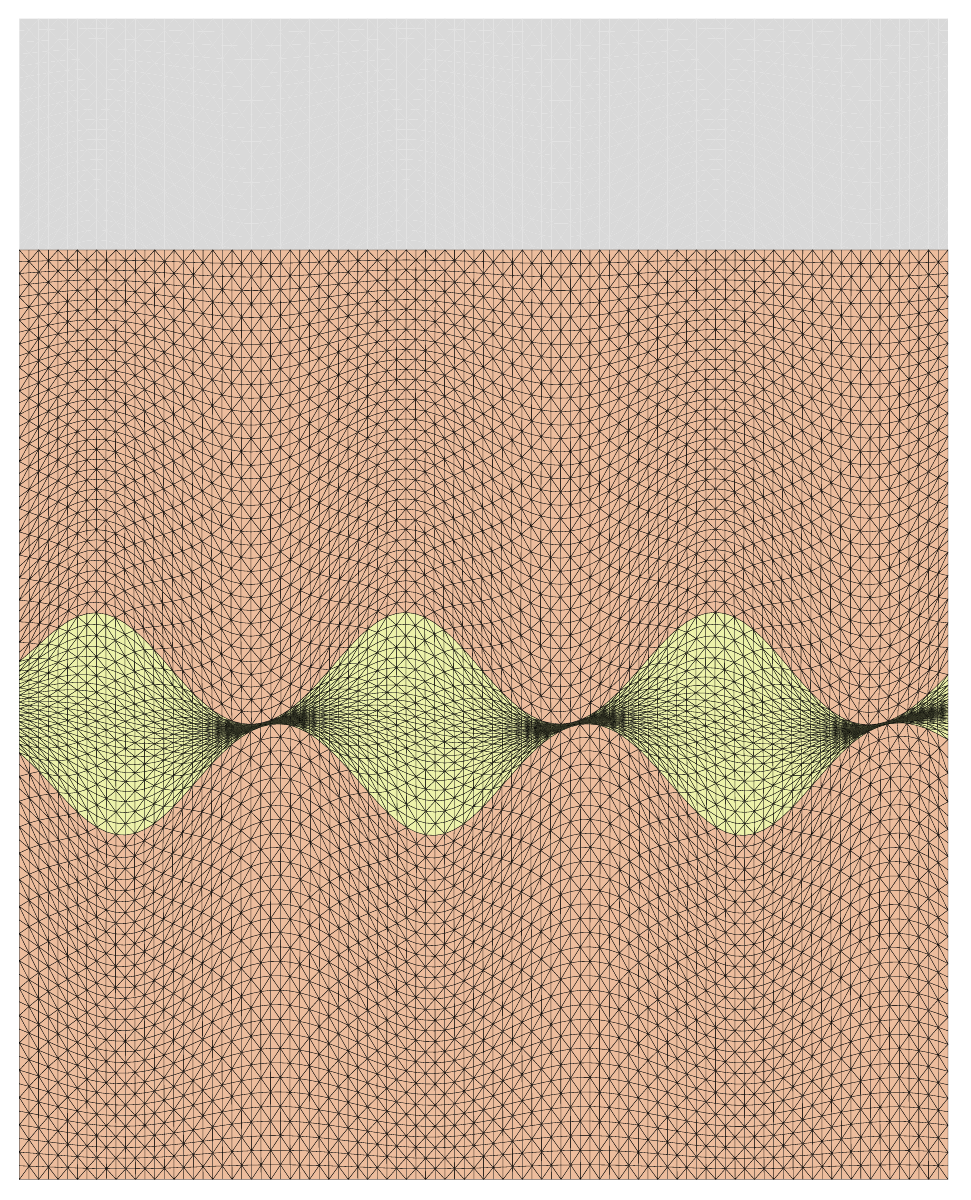}}
\put( 5.7,0.8){\includegraphics[width=4.7cm]{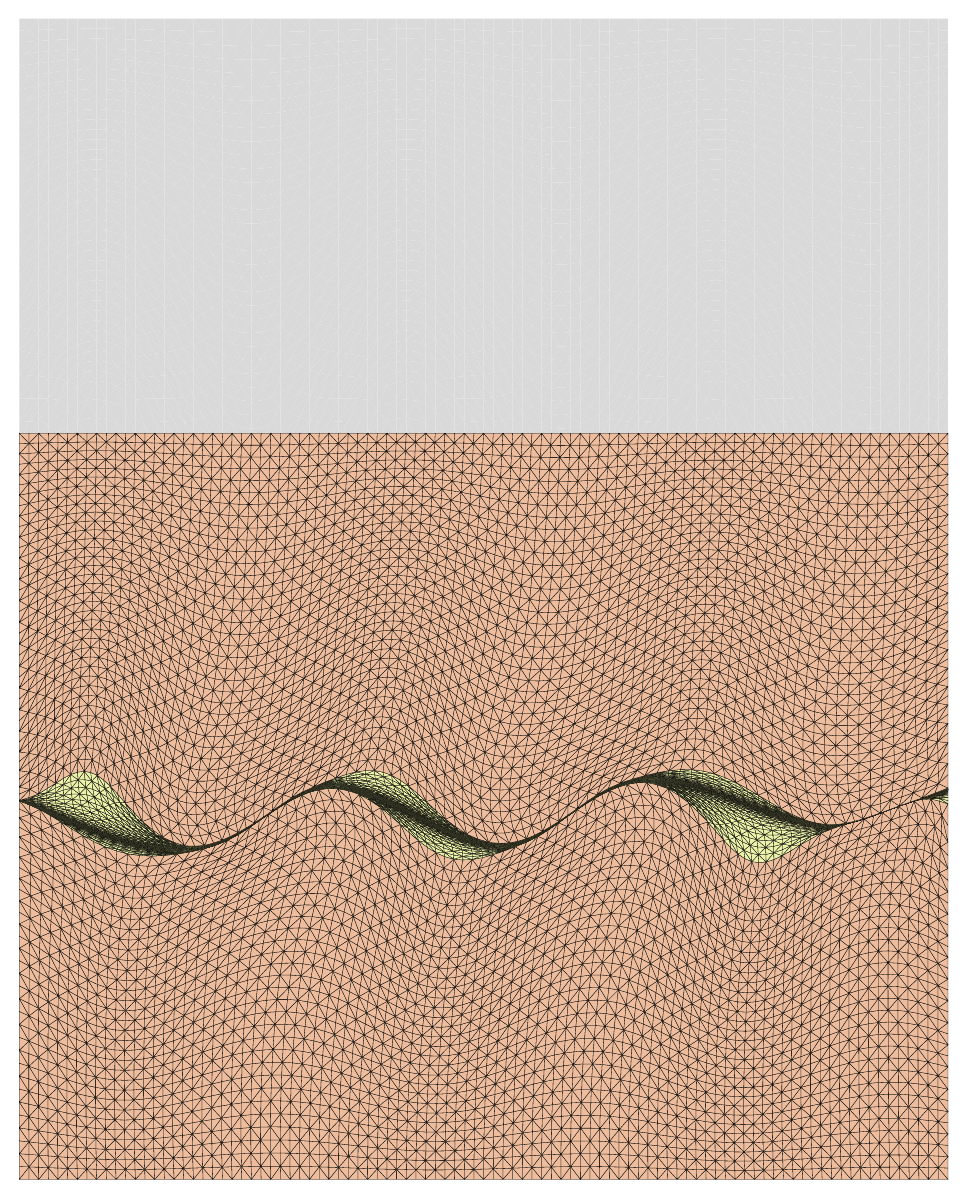}}
\put(10.7,0.8){\includegraphics[width=4.7cm]{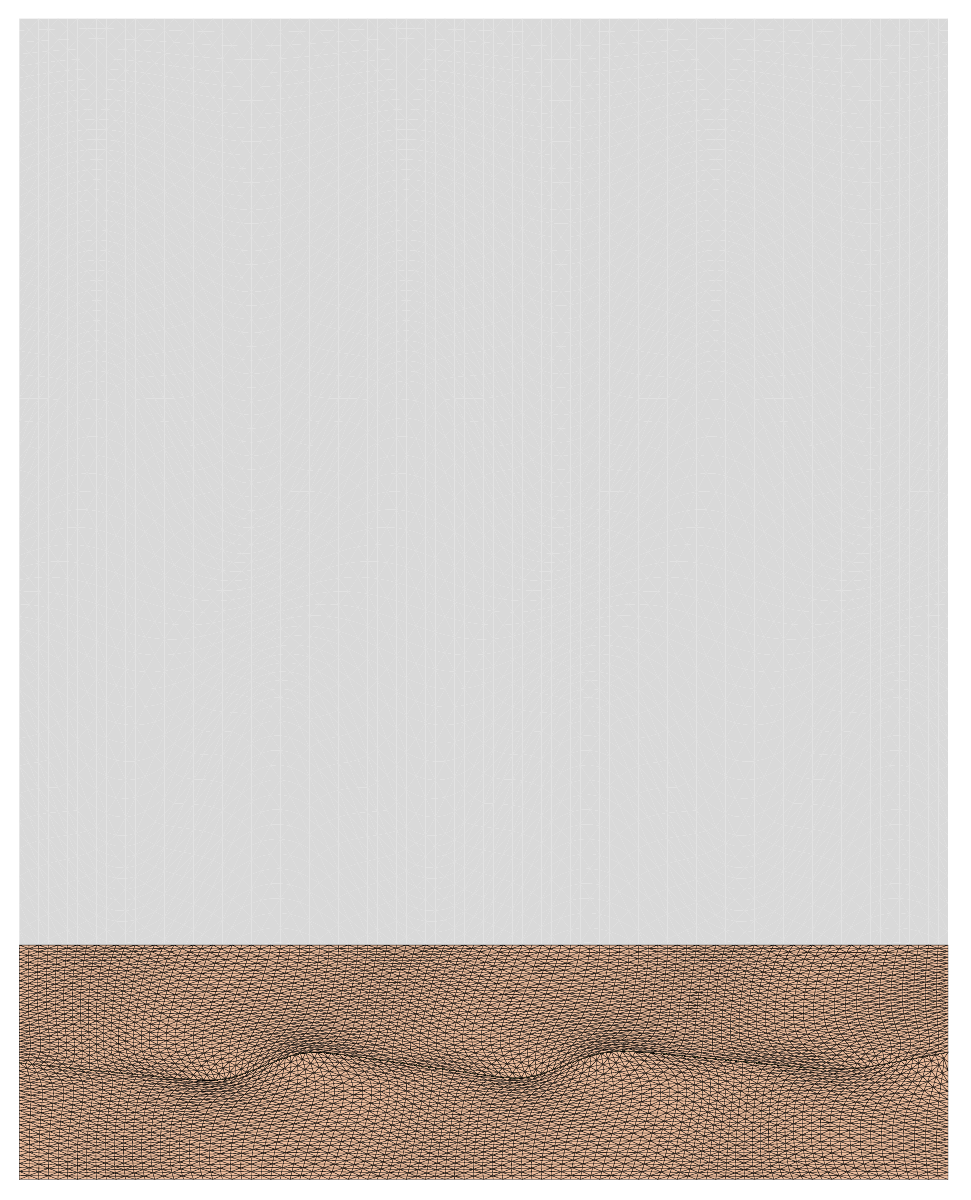}}

\put( 0.4,0.4){a) $\bar u_2=-0.5$mm}
\put( 5.4,0.4){b) $\bar u_2=-0.9$mm}
\put(10.4,0.4){c) $\bar u_2=-2.0$mm}

\end{picture}
\vspace*{-5mm}
\caption{
Deformation of two opposing rough surfaces under increasing compression.
a) shows the initial formation of discrete contact zones,
b) depicts the progressive closure of the intermediate gaps and the growth of the contacting regions,
and c) presents the strongly compressed configuration with an almost continuous contact interface.
The finite element mesh is shown in the deformed configuration.
}
\label{fig:RoughSurfaceCompression}
\end{Figure}

Fig.~\ref{fig:RoughSurfaceCompression} illustrates the progressive closure of the rough interface at three prescribed displacement levels.
For $\bar u_2=-0.5\,\mathrm{mm}$, only isolated asperity pairs are in contact, and the remaining third-medium regions interrupt the direct transport path.
At $\bar u_2=-0.9\,\mathrm{mm}$, the contact zones broaden and additional asperities engage, while the remaining gaps are compressed into narrow localized pockets.
At $\bar u_2=-2.0\,\mathrm{mm}$, the interface is considered as closed and forms a continuous contact path.
The progressive reduction of the third-medium thickness therefore promotes both electrical conduction and heat transfer across the interface.
At the same time, the deformation remains concentrated near the rough surfaces, whereas the surrounding bulk regions deform more uniformly.
The mesh remains connected even under severe local compression, and no visible penetration of the opposing solids occurs.
Consequently, the formulation captures the transition from discrete asperity contact to a fully closed thermo-electrically conducting interface.
Fig.~\ref{fig:RoughSurfaceCompressionConductivity} shows the current transport through the rough interface at initial and full contact.
The contour plot represents the vertical current-density component $J_{\mathrm{e},2}$, while the streamlines indicate the direction of the electric current $\bJ_e$.
At $\bar u_2=-0.5\,\mathrm{mm}$, conduction is restricted to the first contacting tips, which produces strongly localized current paths and pronounced current constriction.
Regions that remain separated by the third-medium do not contribute to the electrical connection.
At $\bar u_2=-2.0\,\mathrm{mm}$, the interface is completely closed and the current crosses the entire contact region.
The current distribution is consequently homogeneous, and the local constriction near isolated contact spots disappears.
These results demonstrate that the proposed switching formulation restricts electrical transport to sufficiently compressed regions of the third medium.

\begin{Figure}[H]
\begin{picture}(0,10.0)
\unitlength1cm

\put( 0.5,0.8){\includegraphics[width=8cm]{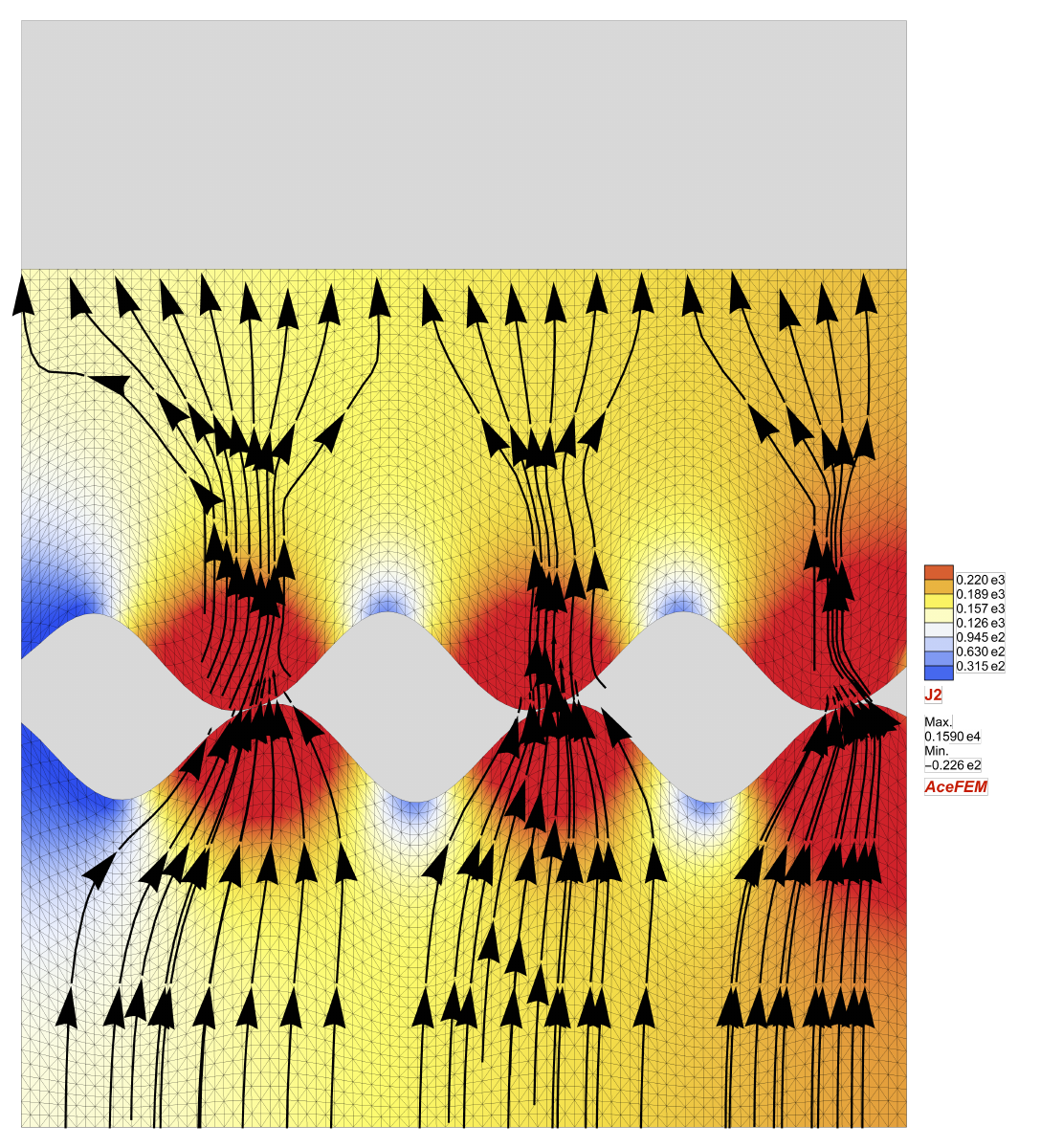}}
\put( 8.0,0.8){\includegraphics[width=8cm]{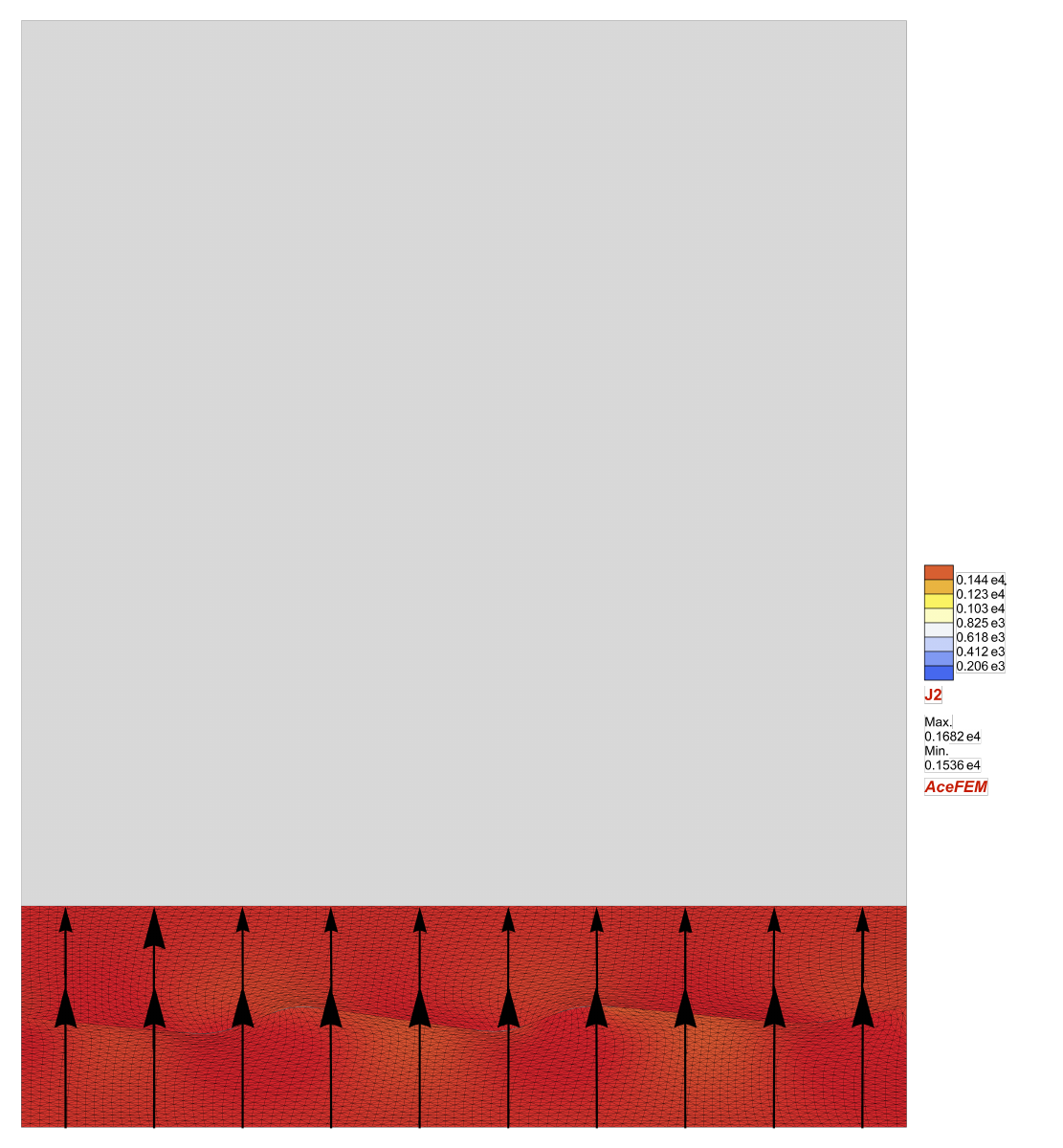}}

\put( 0.5,0.4){a) $\bar u_2=-0.5$mm}
\put( 8.0,0.4){b) $\bar u_2=-2.0$mm}

\end{picture}
\vspace*{-5mm}
\caption{
Current transport through the rough interface at 
a) initial asperity contact and 
b) full compression. 
Contours show the vertical current-density component $J_{\mathrm e,2}$, while
streamlines indicate the current direction. 
}
\label{fig:RoughSurfaceCompressionConductivity}
\end{Figure}

\section{Conclusion and outlook}

A monolithic finite element formulation for electro-thermo-mechanical contact based on the third-medium concept is presented.
Mechanical contact is represented through the deformable intermediate medium without explicit contact detection, gap evaluation or mechanical contact active sets.
Electrical and thermal transport are incorporated through deformation-dependent constitutive switching functions, while Joule heating and the temperature dependence of the electrical conductivity provide the relevant electro-thermal coupling.
Robustness under severe compression is achieved through an independently interpolated deformation-gradient-like field combined with penalty coupling and gradient regularization, which retains compatibility with low-order finite elements.

The numerical examples confirm the expected coupled response.
Electric current remains suppressed before contact and is activated once a sufficiently compressed conducting path forms.
The resulting Joule heat produces a transient temperature increase and, through the temperature-dependent conductivity, a corresponding reduction in current density.
Localized loading leads to strongly non-uniform current and temperature fields, which are resolved without explicit tracking of the evolving contact zone.
For rough interfaces, the formulation captures the transition from isolated asperity contact to an almost continuous conducting interface.
Current constriction is reproduced at the first contact spots, whereas progressive compression distributes the current over an increasing fraction of the nominal contact area.
The examples also demonstrate the decisive influence of the thermal boundary conditions on the predicted temperature evolution.

The present transport laws should be interpreted as deformation-dependent third-medium constitutive models rather than as direct representations of a physical interface gap.
In particular, the electrical activation threshold $J_{\mathrm{crit}}$ and the thermal switching function depend on the chosen third-medium configuration and require problem-specific calibration.
Future work should therefore address their calibration against experimental electrical and thermal contact-resistance data and investigate transport criteria based on geometrically objective interface measures.

\bibliographystyle{unsrtnat}
\let\bibfont\relax
\bibliography{Contact_05}

@Preamble{"\newcommand{\bibfont}{\fontsize{10}{12}\selectfont} "
# "\newcommand{\noopsort}[1]{} "
# "\newcommand{\printfirst}[2]{#1} "
# "\newcommand{\singleletter}[1]{#1} "
# "\newcommand{\switchargs}[2]{#2#1} " }

@STRING{a = {Atherosclerosis}}

@Article{BluSigPou:2021:icm,
  author  = {G. L. Bluhm and O. Sigmund and K. Poulios},
  journal = {Computational Mechanics},
  title   = {Internal contact modeling for finite strain topology optimization},
  year    = {2021},
  pages   = {1099--1114},
  volume  = {67},
  doi     = {10.1007/s00466-021-01974-x},
}

@Article{BluSigPou:2023:ido,
  author  = {G. L. Bluhm and O. Sigmund and K. Poulios},
  journal = {International Journal of Non-Linear Mechanics},
  title   = {Inverse design of mechanical springs with tailored nonlinear elastic response utilizing internal contact},
  year    = {2023},
  pages   = {104552},
  volume  = {157},
  doi     = {10.1016/j.ijnonlinmec.2023.104552},
}

@Article{BogZanKolRan:2015:ncw,
  author  = {T. Bog and N. Zander and S. Kollmannsberger and E. Rank},
  journal = {Computers and Mathematics with Applications},
  title   = {Normal contact with high order finite elements and a fictitious contact material},
  year    = {2015},
  pages   = {1370--1390},
  volume  = {70},
  doi     = {10.1016/j.camwa.2015.04.020},
}

@article{DahSjoDalWal:2026:ara,
  author  = {V. Dahlberg and F. Sj{\"o}vall and A. Dalklint and M. Wallin},
  title   = {A rotation-based approach to third medium contact regularization},
  journal = {Computer Methods in Applied Mechanics and Engineering},
  year    = {2026},
  volume  = {453},
  pages   = {118801},
  doi     = {10.1016/j.cma.2026.118801}
}

@article{DalAleFrePouSig:2025:too,
  author  = {A. Dalklint and J. Alexandersen and A.H. Frederiksen and K. Poulios and O. Sigmund},
  title   = {Topology Optimization of Contact-Aided Thermo-Mechanical Regulators},
  journal = {International Journal for Numerical Methods in Engineering},
  year    = {2025},
  volume  = {126},
  number  = {2},
  pages   = {e7661},
  doi     = {10.1002/nme.7661}
}

@Article{FalAmaHor:2026:dga,
  author  = {O. Faltus and M. Amato and M. Hor{\'a}k},
  journal = {Computer Methods in Applied Mechanics and Engineering},
  title   = {Deformation gradient averaging regularization for third medium contact},
  year    = {2026},
  pages   = {119072},
  volume  = {458},
  doi     = {10.1016/j.cma.2026.119072},
}

@Article{FreDalSigPou:2025:itm,
  author  = {A.H. Frederiksen and A. Dalklint and O. Sigmund and K. Poulios},
  journal = {Computer Methods in Applied Mechanics and Engineering},
  title   = {Improved third medium formulation for 3{D} topology optimization with contact},
  year    = {2025},
  pages   = {117595},
  volume  = {436},
  doi     = {10.1016/j.cma.2024.117595},
}

@article{FreRokPouSigGee:2024:aft,
  author  = {A.H. Frederiksen and O. Roko{\v{s}} and K. Poulios and O. Sigmund and M.G.D. Geers},
  title   = {Adding friction to Third Medium Contact: A crystal plasticity inspired approach},
  journal = {Computer Methods in Applied Mechanics and Engineering},
  year    = {2024},
  volume  = {432},
  pages   = {117412},
  doi     = {10.1016/j.cma.2024.117412}
}

@article{FreSigPou:2024:too,
  author  = {A.H. Frederiksen and O. Sigmund and K. Poulios},
  title   = {Topology optimization of self-contacting structures},
  journal = {Computational Mechanics},
  year    = {2024},
  volume  = {73},
  pages   = {967--981},
  doi     = {10.1007/s00466-023-02396-7}
}

@article{HuaNguZho:2018:aim,
  author  = {J. Huang and N. Nguyen-Thanh and K. Zhou},
  title   = {An isogeometric-meshfree coupling approach for contact problems by using the third medium method},
  journal = {International Journal of Mechanical Sciences},
  year    = {2018},
  volume  = {148},
  pages   = {327--336},
  doi     = {10.1016/j.ijmecsci.2018.08.031}
}

@article{KruNguWriDeL:2018:ifc,
  author  = {R. Kruse and N. Nguyen-Thanh and P. Wriggers and L. De Lorenzis},
  title   = {Isogeometric frictionless contact analysis with the third medium method},
  journal = {Computational Mechanics},
  year    = {2018},
  volume  = {62},
  pages   = {1009--1021},
  doi     = {10.1007/s00466-018-1547-z}
}

@article{WriKorJun:2025:atm,
  author  = {P. Wriggers and J. Korelc and P. Junker},
  title   = {A third medium approach for contact using first and second order finite elements},
  journal = {Computer Methods in Applied Mechanics and Engineering},
  year    = {2025},
  volume  = {436},
  pages   = {117740},
  doi     = {10.1016/j.cma.2025.117740}
}

@article{WriSchSch:2013:afe,
  author  = {P. Wriggers and J. Schr{\"o}der and A. Schwarz},
  title   = {A finite element method for contact using a third medium},
  journal = {Computational Mechanics},
  year    = {2013},
  volume  = {52},
  pages   = {837--847},
  doi     = {10.1007/s00466-013-0848-5}
}

@article{Wri:2026:atm,
  author  = {P. Wriggers},
  title   = {A third medium approach for thermo-mechanical contact based on low order ansatz spaces},
  journal = {Finite Elements in Analysis and Design},
  year    = {2026},
  volume  = {255},
  pages   = {104522},
  doi     = {10.1016/j.finel.2026.104522}
}

@Article{XuWri:2026:sfv,
  author  = {B.-B. Xu and P. Wriggers},
  journal = {Computer Methods in Applied Mechanics and Engineering},
  title   = {Stabilization-free virtual element method for 2{D} third medium contact},
  year    = {2026},
  pages   = {118611},
  volume  = {450},
  doi     = {10.1016/j.cma.2025.118611},
}

@article{XuXueWri:2026:tdt,
  author  = {B.-B. Xu and T. Xue and P. Wriggers},
  title   = {Three-dimensional third medium contact model for hyperelastic contact and pneumatically actuated systems},
  journal = {Journal of the Mechanics and Physics of Solids},
  year    = {2026},
  volume  = {213},
  pages   = {106617},
  doi     = {10.1016/j.jmps.2026.106617}
}

@Article{ZabJanJun:2026:afa,
  author  = {M. von Zabiensky and D.R. Jantos and P. Junker},
  journal = {Finite Elements in Analysis and Design},
  title   = {A fast and robust third medium contact approach using the neighbored element method},
  year    = {2026},
  pages   = {104489},
  volume  = {255},
  doi     = {10.1016/j.finel.2025.104489},
}

@book{Wri:2006:ccm,
  author    = {P. Wriggers},
  title     = {Computational Contact Mechanics},
  publisher = {Springer},
  year      = {2006},
  edition   = {2},
  address   = {Berlin, Heidelberg},
  doi       = {10.1007/978-3-540-32609-0}
}

@book{Lau:2002:cca,
  author    = {T.A. Laursen},
  title     = {Computational Contact and Impact Mechanics: Fundamentals of Modeling Interfacial Phenomena in Nonlinear Finite Element Analysis},
  publisher = {Springer},
  year      = {2002},
  address   = {Berlin, Heidelberg},
  doi       = {10.1007/978-3-662-04864-1}
}

@Book{Yas:2013:nmi,
  author    = {V.A. Yastrebov},
  publisher = {Wiley},
  title     = {Numerical Methods in Contact Mechanics},
  year      = {2013},
  address   = {London},
  isbn      = {9781848215191},
  doi       = {10.1002/9781118647974},
}

@article{PusLau:2004:ams,
  author  = {M.A. Puso and T.A. Laursen},
  title   = {A mortar segment-to-segment contact method for large deformation solid mechanics},
  journal = {Computer Methods in Applied Mechanics and Engineering},
  year    = {2004},
  volume  = {193},
  number  = {6--8},
  pages   = {601--629},
  doi     = {10.1016/j.cma.2003.10.010}
}

@article{SauDeL:2013:acc,
  author  = {R.A. Sauer and L. De Lorenzis},
  title   = {A computational contact formulation based on surface potentials},
  journal = {Computer Methods in Applied Mechanics and Engineering},
  year    = {2013},
  volume  = {253},
  pages   = {369--395},
  doi     = {10.1016/j.cma.2012.09.002}
}

@article{RusTol:2015:dfa,
  author  = {D. Rus and M.T. Tolley},
  title   = {Design, fabrication and control of soft robots},
  journal = {Nature},
  year    = {2015},
  volume  = {521},
  number  = {7553},
  pages   = {467--475},
  doi     = {10.1038/nature14543}
}

@Book{KorWri:2016:aof,
  author    = {J. Korelc and P. Wriggers},
  title     = {{A}utomation of finite element methods},
  year      = {2016},
  publisher = {Springer},
}

@Misc{VorSchWri:2026:arm,
  author        = {M. Vorwerk and J. Schr\"oder and P. Wriggers},
  title         = {A robust mixed finite element formulation for third medium contact},
  year          = {2026},
  archiveprefix = {arXiv},
  eprint        = {2606.28036},
  primaryclass  = {math.NA},
  url           = {https://arxiv.org/abs/2606.28036},
}

@Article{XuXueWri:2026:fov,
  author  = {B.-B. Xu and T. Xue and P. Wriggers},
  journal = {Computer Methods in Applied Mechanics and Engineering},
  title   = {A first-order virtual element method for third-medium contact},
  year    = {2026},
  pages   = {119211},
  volume  = {461},
  doi     = {10.1016/j.cma.2026.119211},
}

@article{WeiWri:2010:nme,
  author  = {C. Wei{\ss}enfels and P. Wriggers},
  title   = {Numerical modeling of electrical contacts},
  journal = {Computational Mechanics},
  volume  = {46},
  pages   = {301--314},
  year    = {2010},
  doi     = {10.1007/s00466-009-0454-8}
}

@article{TerOzhTraShiMatSolKlo:2017:etm,
  author  = {M. Terhorst and O. Ozhoga-Maslovskaja and D. Trauth and A. Shirobokov and P. Mattfeld and M. Solf and F. Klocke},
  title   = {Electro-thermo-mechanical contact model for bulk metal forming under application of electrical resistance heating},
  journal = {The International Journal of Advanced Manufacturing Technology},
  volume  = {89},
  pages   = {3601--3618},
  year    = {2017},
  doi     = {10.1007/s00170-016-9315-8}
}

@article{LiSheKe:2022:mpe,
  author  = {Y.-H. Li and F. Shen and L.-L. Ke},
  title   = {Multi-physics electrical contact analysis considering the electrical resistance and Joule heating},
  journal = {International Journal of Solids and Structures},
  volume  = {256},
  pages   = {111975},
  year    = {2022},
  doi     = {10.1016/j.ijsolstr.2022.111975}
}

@article{LiSheGueKe:2024:aem,
  author  = {Y.-H. Li and F. Shen and M.A. G{\"u}ler and L.-L. Ke},
  title   = {An efficient method for electro-thermo-mechanical coupling effect in electrical contact on rough surfaces},
  journal = {International Journal of Heat and Mass Transfer},
  pages   = {125492},
  year    = {2024},
  doi     = {10.1016/j.ijheatmasstransfer.2024.125492}
}

@Article{WhiCol:1984:hco,
  author  = {G.K. White and S.J. Collocott},
  journal = {Journal of Physical and Chemical Reference Data},
  title   = {Heat Capacity of Reference Materials: Cu and W},
  year    = {1984},
  number  = {4},
  pages   = {1251--1257},
  volume  = {13},
  doi     = {10.1063/1.555728},
}

@Book{TouPowHoKle:1970:tcm,
  author    = {Y.S. Touloukian and R.W. Powell and C.Y. Ho and P.G. Klemens},
  publisher = {IFI/Plenum},
  title     = {Thermal Conductivity: Metallic Elements and Alloys},
  year      = {1970},
  address   = {New York},
  series    = {Thermophysical Properties of Matter: The TPRC Data Series},
  volume    = {1},
}

@Article{Mat:1979:ero,
  author  = {R.A. Matula},
  journal = {Journal of Physical and Chemical Reference Data},
  title   = {Electrical Resistivity of Copper, Gold, Palladium, and Silver},
  year    = {1979},
  number  = {4},
  pages   = {1147--1298},
  volume  = {8},
  doi     = {10.1063/1.555614},
}

@Book{Dav:2001:cac,
  editor    = {J.R. Davis},
  publisher = {ASM International},
  title     = {Copper and Copper Alloys},
  year      = {2001},
  address   = {Materials Park, OH},
  isbn      = {9780871707260},
  series    = {ASM Specialty Handbook},
}


\end{document}